\documentclass[3p]{elsarticle}
\usepackage{hyperref}
\usepackage{amsmath,amssymb,amsthm}
\usepackage{latexsym}            

\usepackage{epsfig}
\usepackage{epstopdf}
\usepackage{amsfonts,amscd, makecell}
\usepackage[ruled,vlined]{algorithm2e}
\usepackage{color,comment}
\usepackage{amsbsy,bm}
\usepackage{subfigure}
\usepackage[utf8]{inputenc}          
\usepackage[T1]{fontenc}

\newcommand\vel{\mathbf{u}}

\newcommand\bx{\boldsymbol{x}}

\newcommand\bn{\mathbf{n}}

\newcommand\bF{\mathbf{F}}

\newcommand\bD{\boldsymbol{D}}

\newcommand\btau{\boldsymbol{\tau}}
\newcommand\bsigma{\boldsymbol{\sigma}}

\newcommand{\be}{\begin{eqnarray}}
\newcommand{\ee}{\end{eqnarray}}
\newcommand{\ben}{\begin{eqnarray*}}
	\newcommand{\een}{\end{eqnarray*}}

\newcommand {\na} {{\nabla}}
\newtheorem{theorem}{Theorem}[section]

\newtheorem{lemma}[theorem]{Lemma}

\newtheorem{rmk}{Remark}[section]
\allowdisplaybreaks \allowdisplaybreaks[4]

\hypersetup{
	colorlinks=true,       
	linkcolor=black,          
	citecolor=green,        
	filecolor=magenta,      
	urlcolor=cyan,         
}

\begin{document}
	\begin{frontmatter}
		\title{An Energy Stable $C^0$ Finite Element Scheme for A Phase-Field Model of Vesicle  Motion  and Deformation}
		
		\author[mymainaddress]{Lingyue Shen}
		\author[mysixthaddress]{Zhiliang Xu}
		\author[mymainaddress]{Ping Lin}	
		\author[mythirdaddress,myfourthaddress,myfifthaddress]{Huaxiong Huang}
		\author[mysecondaryaddress]{Shixin Xu  \corref{mycorrespondingauthor}}
		\cortext[mycorrespondingauthor]{Corresponding author}
		\ead{shixin.xu@dukekunshan.edu.cn}
		\address[mymainaddress]{Department of Mathematics, University of Dundee, Dundee DD1 4HN, Scotland, UnitedKingdom.}
		\address[mysecondaryaddress]{Duke Kunshan University, 8 Duke Ave, Kunshan, Jiangsu, China. }
		\address[mythirdaddress]{Research Centre for Mathematics, Advanced Institute of Natural Sciences, Beijing Normal University (Zhuhai), China;}
		\address[myfourthaddress]{BNU-HKBU United International College, Zhuhai, China}
		\address[myfifthaddress]{Department of Mathematics, University of California, Riverside, 900 University Avenue, Riverside, CA, USA}
		\address[mysixthaddress]{Department of Applied and Computational Mathematics and Statistics, University of Notre Dame,102G Crowley Hall, Notre Dame, IN 46556}
		
		\begin{abstract}	
			A thermodynamically consistent  phase-field model is introduced  for simulating motion and shape transformation of vesicles under flow conditions.
			In particular, a general slip boundary condition is used to describe the interaction between vesicles and the  wall of the fluid domain.
			A second-order accurate  in both space and time  $C^0$ finite element method is proposed to solved the model governing equations. Various
			numerical tests confirm the convergence, energy stability, and conservation of mass and surface area of cells  of the proposed scheme.  Vesicles with different mechanical properties are also used to explain the  pathological risk  for patient with sickle cell disease. 
			
		\end{abstract}
		
	\end{frontmatter}

\section{Introduction}
Studying  dynamic motion and shape transformation of biological cells is always a point of interest in cell  biology, because the shapes of the cells usually relate to their function. For example,  many blood-related diseases are known to be associated with alterations in the geometry
and membrane properties of red blood cells \cite{Takeishi_Ito_2019}. 
Red blood cells in diabetes or sepsis patients exhibit impaired cell deformability \cite{Caimi_Presti_04, Pschl2003EndotoxinBT}. During blood clot formation, an indicator  of platelet activation is its shape change by forming  filopodia and lamellipodia.  Notably, platelets shape changes facilitate their adhesion to the site of vascular injury and cohesion with other platelets or  erythrocytes \cite{shin2017platelet,aslan2012platelet}.

In {\it in silico} study,  it is vitally important to establish a proper  model \cite{lishuwang_2017_vesicle,yangxiaofeng_2015_decoupled, duqiang_2004_bending, Noguchi14159} of cell membranes  for analyzing the dynamical shape transformation of cells in addition to modeling intracellular and extracellular fluids. Various mathematical models  were introduced for predicting cell morphology and function. Dissipative particle dynamics (DPD) \cite{li_2018_dpd} models of red blood cell were developed in \cite{pivkin2009effect,li_2018_dpd,pivkin2008accurate}, and  were used to study effects of red blood cells on platelet aggregation \cite{pivkin2009effect}.
Models based on interface tracking or capturing such as  level set method \cite{xu_level-set_2006,xu_level-set_2014,salac2011level} were also developed \cite{bonito2010parametric,jenkins1977equations,hu2007continuum, hao2015fictitious} to take into consideration the fluid-cell-structure interaction.   In numerical treatment,
 various methods such as immersed boundary method   \cite{laimingzhi_2010_immersed_boundary,laimingzhi_2016_immersed_boundary,laimingzhi_2020_immersed_boundary,wang2020immersed,wu2013simulation}, immersed interface method \cite{hu2016vesicle,kolahdouz2015numerical},   and fictitious domain method \cite{hao2015fictitious} using finite difference or finite element formulation have been introduced to solve governing equations of these models. 


Recently the phase-field approach
has become one of the popular choices for modeling complicated evolution of various structures presented in biological problems \cite{duqiang_2009_variational,yangxiaofeng_2015_decoupled,lishuwang_2017_vesicle} as well as problems in other scientific disciplines \cite{lin2006simulations,Diegel2019finite,wise2009energy,xu2018three}. The phase-field method considers the material interface as a diffuse layer instead of a sharp discontinuity. This regularization can be rigorously formulated through a variational process. The main advantage of the phase-field method is that the phase-field order parameter identifying the diffuse interface is treated as an additional primary unknown of the problem to be solved on the whole domain. Consequently, interface transformations are predicted without the necessity of a remeshing algorithm to treat the evolution of the interface.

Lots of phase-field type vesicle  models have been introduced lately  \cite{wangqi_2011_finite_element,yangxiaofeng_2017_3_component,lowengrub2009phase,du2005phase,du2007analysis,duqiang_2005_curvature,biben2005phase,zhang2009phase,chen2015decoupled,wu2013strong} .  Mechanical properties of the vesicle membrane such as bending stiffness and inextensibility can be modeled rigorously by the phase-field theory \cite{du2005phase,duqiang_2004_bending,duqiang_2008_bending,duqiang_2009_variational} to establish a more comprehensive model. 
For instance, the bending energy $E_{elastic}$ of bending resistance of the lipid bilayer membrane $\Gamma$ in the isotropic case given in the form of the Helfrich bending energy
\be\label{bending energy}
E_{elastic}= \int_{\Gamma} \frac{k}{2} H^2 ds~,
\ee
can be approximated by a modified elastic energy defined on the whole domain in the phase-field formulation \cite{duqiang_2005_curvature,duqiang_2004_bending,duqiang_2008_bending}. 
Here $k$ is the bending modulus and $H$ is the mean curvature of the membrane. Constraints conserving cell mass and ensuring global inextensibility of cell membrane are  frequently introduced into vesicle models to keep the mass and surface area of the vesicle constant \cite{duqiang_2009_variational,voigt_2014_local_inextensibility}. Meanwhile, local inextensibility of vesicle membrane is often needed to prevent from stretching on any point of the vesicle surface \cite{beaucourt_2004_steady}. Aland  {\it et al.} \cite{voigt_2014_local_inextensibility} give out a phase-field based  formulation strongly constraining the local inextensibility of membrane that is defined on the whole domain by introducing a relaxation:
\be\label{voigt reference}
\xi \epsilon^2 \nabla \cdot(\phi^2 \nabla \lambda) +\delta_\epsilon \mathcal{P} : \na\vel=0 & \mbox{in~}\Omega~,
\ee
where $\xi >0$ is a constant controlling the strength of the relaxation, $\phi$ is the phase-field order parameter, $\vel$ is the velocity,  $\delta_\epsilon$ is the delta function which is nonzero on the interface (or membrane) of the vesicle, $\mathcal{P}$ is the tangential projection operator and $\lambda$ can be regarded as a Lagrange multiplier. This  treatment imposes the divergence-free constraint on  $\vel$ only on the interface of the vesicle since in the bulk fluid, the equation degenerates to $\Delta \lambda=0$ and on the interface, it becomes $\mathcal{P} : \na\vel=0$.

The focus of this paper is to model flowing vesicles interacting with the domain boundaries which mimics scenarios such as red blood cells passing a
narrowed blood vessel. This involves considering a moving contact line problem.  
The first goal of this paper thus is to derive a thermodynamically consistent phase-field model for vesicles' motion and shape transformation  in a closed spatial domain by using an energy variational method \cite{shen2020energy,xu2014energetic,xu2018osmosis,guo2015thermodynamically}.  All the physics  taken into consideration are introduced through  definitions of  energy functionals and dissipation functional, together with the kinematic assumptions of laws of conservation.  Besides the energy and dissipation terms defined on bulk region of the domain, terms accounting for boundary effects are also added to the functionals. Then performing variation of these functionals yields an Allen-Cahn-Navier-Stokes (ACNS) system \cite{wu2013strong} with Allen-Cahn general Navier boundary conditions (GNBC) \cite{qian_2006_slipBC}. 
This is in contrast to  most previous works \cite{duqiang_2004_bending, duqiang_2008_bending, chen2015decoupled} in which dynamic boundary condition was  rarely derived during the course of model derivation. Dirichlet or Neumann type conditions were simply added to these models at the end to close the governing equations \cite{voigt_2014_local_inextensibility,duqiang_2009_variational,duqiang_2005_curvature}.
 Moreover, in our model derivation, the incompressibility of the fluid, the local and global inextensibility of the vesicle membrane and the conservation of vesicle mass are taken into account by introducing two Lagrangian multipliers, hydrostatic pressure $p$ and surface pressure $\lambda$ \cite{laimingzhi_2020_immersed_boundary} and penalty terms, respectively.


 The second goal of this paper is to propose an efficient and accurate numerical scheme for solving the obtained fourth-order nonlinear coupled  partial differential equation (PDE) system.   Over the past decades, a lot of schemes have be developed  for Allen-Cahn or Cahn-Hilliard Navier-Stokes systems \cite{yangxiaofeng_2015_decoupled,cheng2019energy,yan2018second,Diegel2019finite,chen2016convergence,guo_2014_midpoint}. Among which,   high-order accurate  energy stable schemes such as the Invariant Energy Quadratization (IEQ) Schemes \cite{yang2017numerical,yang2020convergence} and the scalar auxiliary variable (SAV)  method \cite{shen2018convergence,shen2018scalar,shen2019new,gong2019arbitrarily,li2020sav} are proposed based on approximate energy expression. 
  As for systems such as vesicle models introduced in the current and other works which are more sophisticated than Allen-Cahn or Cahn-Hilliard Navier-Stokes systems, backward Euler time discretization method is frequently used \cite{voigt_2014_local_inextensibility,duqiang_2004_bending,gu2016two,gu2014simulating} leading to a first-order accurate scheme. Later on, decoupled energy stable schemes  are proposal by Chen $\&$ Yang in \cite{chen2015decoupled}, and   Francisco $\&$ Giordano \cite{guillen2018unconditionally} by introducing explicit, convective velocities.  
In the current work,  an efficient, energy-law preserving (thus energy stable) and  second-order accurate $C^0$ finite element  (FE) scheme is proposed to solve the obtained vesicle system using ideas introduced in \cite{guo_2014_midpoint}. The key idea of this scheme is to utilize the mid-point method in time discretizaiton to ensure the accuracy in time, and the form of the law of the discrete energy dissipation is  same as that of the continuous model. In order to properly treat the term related to inextensibility of the membrane,  a relaxation term of local inextensibility as in \cite{voigt_2014_local_inextensibility} is introduced.   The numerical study of convergence confirms the proposed scheme is second-order convergence in both time and space. Furthermore vesicle deformation simulations illustrate it is  
energy stable, and numerically conserves mass and surface area of vesicles. 


The introduction of the GNBC in this work makes it possible to study the more complicated fluid-structure interaction problems. In this paper, the developed model is applied to studying vesicles passing narrow channels. The results confirm  that the more rounded the vesicles (smaller surface-volume ratio) are, the more likely the vesicles form lockage when they pass through narrow channels. 
 It is also worth noting that it is critical to include the local inextensibility of the vesicle membrane in the model when studying this type of problems. Without the local inextensibility, the vesicle membrane can be falsely stretched or compressed. 
Lastly, although membrane structures of vesicles and blood cells are quite different, a blood cell in many studies can be treated as an elastic capsule with bending rigidity, in which the membrane is impenetrable to both interior and exterior fluids. Therefore our model developed for vesicles can be readily applied for studying a vast body of blood cells related problems \cite{marth2016margination}.  

The rest of paper is organized as follows. 
Section 2 of the paper  begins with introducing basic dynamical assumptions that have been used in many papers \cite{duqiang_2009_variational,qian_variational_2006}, and is devoted to model derivation. Dimensionless model governing equations and the  energy decaying law of the model are presented in Section 3. In Section 4,  the numerical scheme solving the proposed model is developed, and  its energy law is given. Numerical simulation results are described in Section 5 to confirm the  energy law of the numerical scheme and the feasibility of our model. A case study of vesicle passing through a narrow channel is shown, which is to simulate the motion of red blood cells in blood vessel. Conclusions are drawn in Section 6.


\section{Model Derivation}
Derivation of the model for simulating a flowing vesicle deforming in a channel filled with extracellular fluid is presented in this section.
The phase-field label function $\phi$  is used to track the location of the cell membrane, where $\phi(\boldsymbol x)=1$ denotes the intracellular space and $\phi(\boldsymbol x) =-1$ denotes the extracellular space of the vesicle. In the diffuse interface layer, $\phi(\boldsymbol x) \in (-1, 1)$.  The derived model obeys the following energy dissipation law $\frac{dE_{total}}{dt} = -\Delta$ \cite{xu2014energetic,eisenberg2010energy}, where $E_{total}$ is the total energy of the system being modeled, and $\Delta$ is the energy dissipation (rate). This energy law describes how  energy of a system dissipates when it undergoes an isothermal process, and can be easily derived from the first and second laws of the thermodynamics. 

Model kinematic assumptions in the channel domain $\Omega$ are given as follows:
\be\left\{\begin{array}{l}
\label{assumption}
\frac{\partial \phi}{\partial t} +\nabla\cdot (\vel \phi) = q_{\phi}~,\\
\rho(\frac{\partial\vel}{\partial t}+(\vel\cdot\nabla)\vel) = \nabla\cdot\bsigma_{\eta}+\bF_{\phi}~,\\
\nabla \cdot\vel = 0~,\\
 \delta  \mathcal{P} : \na\vel=0~, \end{array}\right.
\ee
where $\vel$ is the macroscale  fluid velocity. $q_{\phi}$ which eventually takes form of generalized diffusion is a phase-field phenomenological term,  $\bsigma_{\eta}$ is fluid shear stress, and $\bF_{\phi}$ is the cell membrane force acting on the fluid. The projection operator $\mathcal{P}$ is defined to be $(I-\bn_{m}\otimes \bn_{m})$, where $\bn_m=\frac{\nabla\phi}{|\nabla\phi|}$ is the unit outward normal vector of the  membrane or interface.  
$q_{\phi}$, $\bsigma_{\eta}$ and $\bF_{\phi}$ are yet to be determined. 

In the system \eqref{assumption},  the first equation is the Allen-Cahn type of phase-field equation to track the interface. 
The second equation is  the conservation of momentum. The third equation accounts for the fluid incompressibility. 
The last  equation  is a relaxation of local inextensibility  condition $\nabla_{\Gamma}\cdot \vel = 0$ on the interface $\Gamma$ \cite{lowengrub2009phase,marth2016margination}.  $\delta$ is the regularized dirac delta function supported on the diffuse interface layer
\cite{lee2012regularized}.

On the boundary $\partial\Omega_w$ of the domain,  the following boundary conditions are assumed
\be\left\{\begin{array}{l}\label{assumption_bd}
	\vel\cdot\bn = 0~, \\ 
	\vel_{\tau}\cdot \btau_i = f_{\tau_i}~,\\
	\frac{\partial \phi}{\partial t} + \vel\cdot\nabla_{\Gamma}\phi = J_{\Gamma}~, \\
	f = 0~,\\
	\partial_n\lambda = 0~,
\end{array}\right.\ee
where  an Allen-Cahn type boundary condition is employed for $\phi$, $\vel_{\tau} = \vel-(\vel\cdot\bn)\bn$ is the fluid slip velocity with respect to the wall, $\btau_i, i=1,2$ are the tangential directions on the wall surface (2D), and $\nabla_{\Gamma} = \nabla-\bn(\bn\cdot\nabla)$ is the surface gradient operator on the boundary $\partial\Omega_w$. 

Rest of this section is devoted to deriving the exact forms of $ q_{\phi}$, $\bsigma_{\eta}$, $\bF_{\phi}$ and $J_{\Gamma}$ using the energy variational method \cite{shen2020energy,xu2018osmosis}.
By following the works in \cite{wu2013strong,du2005phase}, the total energy functional of a cell- (or vesicle-) fluid system is defined as follows:
\be\label{energy}
E_{total}&=& \underbrace{E_{kin}}_{Macroscale}+\underbrace{E_{cell}+E_{w}}_{Microscale}~,
\ee
where $E_{kin}$ is the kinetic energy, $E_{cell}$ is the cell membrane energy and
$E_w$ is the specific wall energy due to cell-wall interaction. 

The kinetic energy accounts for the transport of the  cell-fluid mixture, and is defined as:
\be
E_{kin}=\int_{\Omega}\left(\frac 1 2 \rho | \vel|^2\right)d \boldsymbol{x}~,
\ee
where  $\rho$ is the macroscale density of the mixture, and is assumed to be equal to a constant $\rho_0$ in this work (matched density case). 

If  the cell membrane is assumed to be iostropic and is only composed of lipid bilayer, the bending energy of the bending resistance of the cell membrane can be modeled by the Helfrich bending energy \cite{du2005phase}. In the phase-field formulation, the Helfrich bending energy can be regularized as follows
\be
E_{bend} = \int_{\Omega} \frac{\hat\kappa_B}{2\gamma}\left|\frac{f(\phi)}{\gamma}\right|^2d\bx~,
\ee
where  $\hat\kappa_B$ is the bending modulus, $f(\phi)=\frac{\delta G}{\delta\phi} = -\gamma^2 \Delta\phi+ (\phi^2-1)\phi$, $G(\phi)= \frac{\gamma^2|\nabla\phi|^2}2+\frac{(1-\phi^2)^2}{4}$, and $\gamma$ is related to the membrane thickness.  

Penalty terms $\frac{M_v}{2}\frac{\left(V(\phi)-V(\phi_0)\right)^2}{V(\phi_0)}$ and $\frac{M_s}{2}\frac{\left(S(\phi)-S(\phi_0)\right)^2}{S(\phi_0)}$ are introduced into the membrane energy functional 
in order to preserve the total volume and surface area of the cell, respectively. To this end, the cell energy $E_{cell}$ is defined to be
\be
E_{cell}
&=&E_{bend}+\frac{M_v}{2}\frac{\left(V(\phi)-V(\phi_0)\right)^2}{V(\phi_0)} +\frac{M_s}{2}\frac{\left(S(\phi)-S(\phi_0)\right)^2}{S(\phi_0)}~,
\ee
where $V(\phi) = \int_{\Omega}\phi d \boldsymbol x$ is the volume difference of the cell-fluid system and $S(\phi) = \int_{\Omega} \frac{G(\phi)}{\gamma} d\boldsymbol x$ is the surface area of the cell.
$M_v$ and $M_s$ are   cell volume and surface area constraint coefficients, respectively.

The wall free energy $E_w$   is defined on the channel wall $\partial\Omega_w$ to account for the interaction between the cell, the fluid and the wall, and adopts the following form \cite{qian_molecular_2006,qian_variational_2006}
\be
E_w=\sigma \int_{\partial\Omega_w}f_w(\phi)ds ~, 
\ee
where  $f_w$ is given by
\be
f_w(\phi)= -\sigma \sin\left(\frac{\phi\pi}{2}\right)\cos(\theta_s)~.
\ee
Here $\sigma$ is the surface tension of cell  and  $\theta_s$   is the static contact angle \cite{ren2007boundary,ren_heterogeneous_2005}.

The chemical potential $\mu$ is obtained by taking the variation of $E_{bulk} = E_{kin}+E_{cell}$ with respect to $\phi$,
\be\label{mu}
\mu = \frac{\delta E_{bulk}}{\delta \phi} = \frac{\hat\kappa_B}{\gamma^3}g(\phi) +M_v\frac{V(\phi)-V(\phi_0)}{V(\phi_0)}+\frac{M_s}{\gamma}\frac{S(\phi)-S(\phi_0)}{S(\phi_0)}f(\phi)~, 
\ee
where $g(\phi) = -\gamma^2\Delta f+(3\phi^2-1)f(\phi)$. 

  It is assumed in the present work that  dissipation of the system energy is due to fluid viscosity, friction on the wall, and interfacial mixing due to diffuse interface representation.    Accordingly, the total dissipation functional $\Delta$ is defined as follows
\be\label{dispp}
\Delta &=& \int_{\Omega}2\eta|\bD_{\eta}|^2d\boldsymbol x+\int_{\Omega}\frac 1 M_{\phi}|q_{\phi}|^2d\boldsymbol x+  \int_{\partial\Omega_w}\beta_s|\vel_{\btau}|^2ds +\int_{\partial\Omega_w}\kappa_{\Gamma}|J_{\Gamma}|^2ds~, 
\ee
 where the first term is the macroscopic dissipation induced by the fluid viscosity, the second term is the microscopic dissipation induced by the diffuse interface,  the third term is the boundary friction dissipation and the last term is the dissipation induced by   the diffuse interface   contacting  the wall.

By taking the time derivative of the total energy functional \eqref{energy}, it is  obtained that
\be\label{dEdt}
\frac{d E_{total}}{dt} &=& \frac{d}{dt} E_{kin} + \frac{d}{dt} E_{cell} + \frac{d}{dt} E_w\\
& \equiv & I_1+I_2+I_3~.
\ee

This leads to
\be\label{I_1}
I_1 &=&\frac{d}{dt}\int_{\Omega}\frac{\rho|\vel|^2}{2}d\bx \nonumber\\
&=&\int_{\Omega}\frac 1 2 \frac{\partial\rho}{\partial t}|\vel|^2d\bx +\int_{\Omega} \rho\frac{\partial \vel}{\partial t}\cdot\vel d\bx\nonumber\\
&=&\int_{\Omega}\frac 1 2 \frac{\partial\rho}{\partial t}|\vel|^2d\bx +\int_{\Omega} \rho\frac{d \vel}{d t}\cdot\vel d\bx-\int_{\Omega} \left(\rho \vel\cdot\nabla \vel\right)\cdot \vel d\bx\nonumber\\
&=&\int_{\Omega}\frac 1 2 \frac{\partial\rho}{\partial t}|\vel|^2d\bx +\int_{\Omega} \rho\frac{d \vel}{d t}\cdot\vel d\bx+\int_{\Omega} \nabla\cdot(\rho\vel)\frac{|\vel|^2}2d\bx \nonumber\\
&=&\int_{\Omega}( \nabla \cdot \bsigma_\eta)  \cdot\vel d\bx +\int_{\Omega}\bF_{\phi}\cdot \vel d\bx +\int_\Omega \lambda\delta \mathcal{P} : \na\vel d\bx -\int_{\Omega} pI:\nabla \vel d\bx \nonumber\\
&=&-\int_{\Omega}((\bsigma_{\eta}+pI):\nabla\vel)d\bx+\int_{\Omega}\bF_{\phi}\cdot \vel dx -\int_\Omega \na\cdot(\lambda\delta  \mathcal{P})\cdot\vel d\bx \nonumber\\
&&+\int_{\partial\Omega_w}((\bsigma_{\eta}+\lambda\delta \mathcal{P})\cdot\bn)\cdot\vel_{\tau} dS~,
\ee
where $p$ and  $\lambda$ are  introduced as Lagrange multipliers accounting for fluid incompressiblity and local inextensibility of the cell membrane, respectively.  
The velocity boundary conditions in \eqref{assumption_bd} and integration by parts are utilized in the above derivation from step 4 to  step 6.

Using the first equation in Eq.~set~\eqref{assumption},  and the definitions of $g(\phi)$ and $f(\phi)$ give rise to 
\be\label{I_2}
I_2 &=& \frac{d}{dt}\int_{\Omega}\frac{\hat\kappa_B}{2\gamma}\left|\frac{f(\phi)}{\gamma}\right|^2d\bx +\frac{d}{dt}\left(\frac{M_v}{2}\frac{\left(V(\phi)-V(\phi_0)\right)^2}{V(\phi_0)} +\frac{M_s}{2}\frac{\left(S(\phi)-S(\phi_0)\right)^2}{S(\phi_0)}\right)\nonumber\\
&=& \int_{\Omega}\frac{\hat\kappa_B}{\gamma}\frac{f}{\gamma^2}\frac{\partial f}{\partial t}d\bx +M_v\int_{\Omega}\frac{V(\phi)-V(\phi_0)}{V(\phi_0)}\frac{\partial\phi}{\partial t}d\bx +M_s\int_{\Omega}\frac{S(\phi)-S(\phi_0)}{S(\phi_0)}\frac{\partial S(\phi)}{\partial t}d\bx\nonumber\\
&=& \int_{\Omega}\frac{\hat\kappa_B}{\gamma}\frac{f}{\gamma^2}\left(-\gamma^2\Delta\left(\frac{\partial \phi}{\partial t}\right) +(3\phi^2-1)\frac{\partial \phi}{\partial t}\right)d\bx\nonumber\\ &&+M_v\int_{\Omega}\frac{V(\phi)-V(\phi_0)}{V(\phi_0)}\frac{\partial\phi}{\partial t}d\bx \nonumber\\
&&+M_s\int_{\Omega}\frac{S(\phi)-S(\phi_0)}{S(\phi_0)}\frac{1}{\gamma}\left(\gamma^2\nabla\phi\cdot \nabla\frac{\partial \phi}{\partial t} +(\phi^2-1)\phi \frac{\partial \phi}{\partial t} \right)d\bx\nonumber\\
&=& \int_{\Omega}\frac{\hat\kappa_B}{\gamma^3}\left(-\gamma^2\Delta f +(3\phi^2-1)f\right)\frac{\partial \phi}{\partial t}d\bx\nonumber\\ &&+M_v\int_{\Omega}\frac{V(\phi)-V(\phi_0)}{V(\phi_0)}\frac{\partial\phi}{\partial t}d\bx \nonumber\\
&&+M_s\int_{\Omega}\frac{S(\phi)-S(\phi_0)}{S(\phi_0)}\frac{1}{\gamma}\left(-\gamma^2\Delta\phi +(\phi^2-1)\phi   \right)\frac{\partial \phi}{\partial t}d\bx\nonumber\\
&&-\int_{\partial\Omega_w}\frac{\hat\kappa_B}{\gamma}f\frac{\partial}{\partial t}(\partial_n\phi) ds +\int_{\partial\Omega_w} \frac{\hat\kappa_B}{\gamma} \partial_n f\frac{\partial \phi}{\partial t}ds +M_s\int_{\partial\Omega_w} \frac{S(\phi)-S(\phi_0)}{S(\phi_0)}\gamma\partial_n\phi\frac{\partial \phi}{\partial t}ds\nonumber\\
&=&\int_{\Omega}\mu\frac{\partial \phi}{\partial t} d\bx -\int_{\partial\Omega_w}\frac{\hat\kappa_B}{\gamma}f\frac{\partial}{\partial t}(\partial_n\phi) ds\nonumber\\
&&+\int_{\partial\Omega_w} \frac{\hat\kappa_B}{\gamma} \partial_n f\frac{\partial \phi}{\partial t}ds +M_s\int_{\partial\Omega_w} \frac{S(\phi)-S(\phi_0)}{S(\phi_0)}\gamma\partial_n\phi\frac{\partial \phi}{\partial t}ds\nonumber\\
&=& \int_{\Omega}\mu q_{\phi}d\bx-\int_{\Omega}\mu\vel\cdot\nabla \phi d\bx\nonumber\\
 &&+\int_{\partial\Omega_w} \frac{\hat\kappa_B}{\gamma} \partial_n f\frac{\partial \phi}{\partial t}ds +M_s\int_{\partial\Omega_w} \frac{S(\phi)-S(\phi_0)}{S(\phi_0)}\gamma\partial_n\phi\frac{\partial \phi}{\partial t}ds~.
\ee 
Here the second and third equations in the boundary conditions \eqref{assumption_bd} and integration by parts are used in step 4 in the above derivation. 

For $I_3$ in Eq.~(\ref{dEdt}), it is easy to see that 
\be\label{I_3}
I_3 = \int_{\partial\Omega_w}\frac{\partial f_w}{\partial \phi}\frac{\partial \phi}{\partial t}ds~. 
\ee

Combining Eqs. \eqref{I_1} to \eqref{I_3} yields
\be
\frac{d}{dt}E_{total}
&=&-\int_{\Omega}((\bsigma_{\eta}+pI):\nabla\vel)d\bx+\int_{\Omega}(\bF_{\phi}-\mu\nabla\phi-\na\cdot(\lambda\delta \mathcal{P}))\cdot \vel d\bx+\int_{\Omega}\mu q_{\phi}d\bx \nonumber\\
&&+\int_{\partial\Omega_w}((\bsigma_{\eta}+\lambda\delta  \mathcal{P})\cdot\bn)\cdot\vel_{\tau} ds  +\int_{\partial\Omega_s}\hat L(\phi)\frac{\partial \phi}{\partial t} ds \nonumber\\
&=&-\int_{\Omega}((\bsigma_{\eta}+pI):\nabla\vel)d\bx+\int_{\Omega}(\bF_{\phi}-\mu\nabla\phi-\na\cdot(\lambda\delta_\epsilon \mathcal{P}))\cdot \vel d\bx+\int_{\Omega}\mu q_{\phi}d\bx \nonumber\\
&&+\int_{\partial\Omega_w}((\bsigma_{\eta}+\lambda\delta(\phi) \mathcal{P})\cdot \bn )\cdot\vel_{\tau} ds+\int_{\partial\Omega_s}\hat L(\phi)(-\vel\cdot\nabla_{\Gamma}\phi+J_{\Gamma})ds\nonumber\\
&=&-\int_{\Omega}((\bsigma_{\eta}+pI):\nabla\vel)d\bx+\int_{\Omega}(\bF_{\phi}-\mu\nabla\phi-\na\cdot(\lambda\delta  \mathcal{P}))\cdot \vel d\bx+\int_{\Omega}\mu q_{\phi}d\bx \nonumber\\
&&+\int_{\partial\Omega_w}((\bsigma_{\eta}+\lambda\delta  \mathcal{P})\cdot\bn-\hat L(\phi)\nabla_{\Gamma}\phi)\cdot\vel_{\tau} ds+\int_{\partial\Omega_w}\hat L(\phi)J_{\Gamma}ds~,\nonumber\\
\ee
where $ \hat L(\phi) = \displaystyle\frac{\hat\kappa_B}{\gamma}\partial_nf +M_s\frac{S(\phi)-S(\phi_0)}{S(\phi_0)}\gamma\partial_n\phi+\frac{\partial f_w}{\partial \phi}$.

 Using the energy dissipation law $\frac{dE_{total}}{dt} = -\Delta$ \cite{xu2014energetic,eisenberg2010energy}, and the definition of the dissipation functional \eqref{dispp}, it is obtained that 
\be\left\{\begin{array}{ll}
\bsigma_{\eta} = 2\eta\bD_{\eta}-pI~,& \mbox{in~} \Omega~,\\
q_{\phi} = -M_{\phi}\mu~, &\mbox{in~} \Omega~,\\
\bF_{\phi} = \mu\nabla\phi+\nabla\cdot(\lambda\delta \mathcal{P})~, & \mbox{in~}\Omega~,\\
J_{\Gamma} = -\kappa_{\Gamma}^{-1} \hat L(\phi)~, &\mbox {on~}\partial\Omega_w~,\\
u_{\tau_i} =\beta_{s}^{-1}(-(\bn\cdot(\bsigma_{\eta}+\lambda\delta_\epsilon \mathcal{P})\cdot\btau_i)+\hat L(\phi)\partial_{\tau_i}\phi)~,~ i=1,2,& \mbox{on~} \partial\Omega_w~.
\end{array}\right.
\ee
 Here constant $M_\phi$ is called the mobility (a phenomenological parameter), $\kappa_{\gamma}$ is the boundary mobility (a phenomenological parameter) and 
 $\beta_s$ is the wall friction coefficient.


To this end,  the proposed phase-field model is composed of  the following equations 
\be\left\{\begin{array}{l}
\label{sys}
\frac{\partial \phi}{\partial t} +\nabla\cdot (\vel \phi) = -M_{\phi}\mu~,\\
\mu = \frac{\hat\kappa_B}{\gamma^3}g(\phi) +M_v\frac{V(\phi)-V(\phi_0)}{V(\phi_0)}+\frac{M_s}{\gamma}\frac{S(\phi)-S(\phi_0)}{S(\phi_0)}f(\phi)~,\\
g(\phi) = -\gamma^2\Delta f+(3\phi^2-1)f(\phi),\\ f(\phi)= -\gamma^2 \Delta\phi+ (\phi^2-1)\phi~,\\
\rho(\frac{\partial\vel}{\partial t}+(\vel\cdot\nabla)\vel)+\nabla p = \nabla\cdot(2\eta\bD_{\eta})+\mu\nabla\phi+\na\cdot(\lambda\delta   \mathcal{P})~,\\
\nabla \cdot\vel = 0~,\\
\delta  \mathcal{P} : \na\vel=0~,
\end{array}\right.
\ee
with the boundary conditions
\be\left\{\begin{array}{l}
\label{sys_bd}
\vel\cdot\bn =0~,\\
-\beta_{s}u_{\tau_i} =(\bn\cdot(\bsigma_{\eta}+\lambda\delta \mathcal{P})\cdot\btau_i)-\hat   L(\phi)\partial_{\tau_i}\phi~, ~i=1,2,\\
f=0~,\\
\kappa_{\Gamma}\left(\frac{\partial \phi}{\partial t} + \vel\cdot\nabla_{\Gamma}\phi\right) = -\hat L(\phi)~,\\
\hat L(\phi) = \frac{\hat\kappa_B}{\gamma}\partial_nf +M_s\frac{S(\phi)-S(\phi_0)}{S(\phi_0)}\gamma\partial_n\phi+\frac{\partial f_w}{\partial \phi}~,\\
\partial_n\lambda = 0~.
 \end{array}
\right.
\ee

\section{Dimentionless Model Governing Equations and Energy Dissipation Law}

If the  viscosity, length, velocity, time, bulk and boundary chemical potentials in Eqs.~(\ref{sys})-(\ref{sys_bd}) are scaled by their corresponding characteristic values $\eta_0$, $L$, $U$, $\frac{L}{U}$ $\frac{\eta_0 U}{L}$ and $\eta_0U$, respectively,  Eqs.~(\ref{sys})-(\ref{sys_bd}) can be rewritten as 

\be\left\{\begin{array}{ll}
\label{sys_nd}
	Re(\frac{\partial \vel}{\partial t}+(\vel\cdot \na)\vel)+\na P=\na\cdot(2\eta \bD)+\mu\na\phi+\na\cdot(\lambda\delta_\epsilon \mathcal{P})~,& \mbox{in~}\Omega~,\\[3mm]
	\na\cdot \vel=0~,& \mbox{in~}\Omega~,\\[3mm]
	\frac{\partial\phi}{\partial t}+\vel\cdot\na\phi=-\mathcal{M}\mu~,& \mbox{in~}\Omega~,\\[3mm]
	\mu=     \kappa_B g(\phi) +\mathcal{M}_v\frac{(V(\phi)-V(\phi_0))}{V(\phi_0)}+\mathcal{M}_s\frac{(S(\phi)-S(\phi_0))}{S(\phi_0)}f(\phi)~,  & \mbox{in~}\Omega~,\\[3mm]
	f(\phi) = -\epsilon \Delta\phi+ \frac{(\phi^2-1)}{\epsilon}\phi,~g(\phi) = -\Delta f+\frac{1}{\epsilon^2}(3\phi^2-1)f(\phi)~, & \mbox{in~}\Omega~,\\[3mm]
\delta_\epsilon \mathcal{P} : \na\vel=0 ~, & \mbox{in~}\Omega~,
\end{array}\right.
\ee
with  the boundary conditions 
\be
\displaystyle \left\{\begin{array}{ll}
\label{bd_nd}
	\kappa \dot{\phi}+L(\phi)=0~, & \mbox{on~} \partial\Omega_w~,\\
	L(\phi)= \kappa_B \partial_n f+ \epsilon\mathcal{M}_s\frac{S(\phi)-S(\phi_0)}{S(\phi_0)} \partial_n\phi+ \alpha_w\frac{df_w}{d\phi}~,  & \mbox{on~} \partial\Omega_w~,\\
	-l_s^{-1} u_{\tau_i} = \boldsymbol{\tau_i}\cdot( 2\eta \mathbf{D}_{\eta}+\lambda\delta_\epsilon \mathcal{P})\cdot \boldsymbol{n} - L(\phi)\partial_{\tau_i}\phi~,~i=1,2, & \mbox{on~}  \partial\Omega_w~,\\
  f=0~, &\mbox{on~} \partial\Omega_w~,\\
  \partial_n\lambda = 0~, &\mbox{on~} \partial\Omega_w~,
\end{array}\right.
\ee
where $V(\phi) = \displaystyle\int_{\Omega}\phi d\bx$ and $S(\phi) = \displaystyle\int_{\Omega}\frac{\epsilon}{2}|\nabla\phi|^2 +\frac{1}{4\epsilon}(\phi^2-1)^2 d\bx$. 
The dimensionless constants appeared in Eqs.~(\ref{sys_nd})-(\ref{bd_nd}) are given by $\epsilon = \frac{\gamma}{L}$, $Re = \frac{\rho_0 UL}{\eta_0}$, $\mathcal{M} = M_\phi\eta_0$, 
$\kappa_B = \frac{\hat\kappa_B}{L^2\eta_0U}$, 
$k = \frac{\hat\kappa_B  }{\eta_0 L}$, 
$l_s = \frac{\eta_0}{\beta_{s}L}$, $\alpha_w = \frac{\sigma}{\eta_0 U}$, $\mathcal{M}_s = \frac{M_s}{\eta_0 U}$, and $\mathcal{M}_v = \frac{M_vL}{\eta_0 U}$. 

Let  $\| f \| = \left(\int_{\Omega}|f|^2 d\bx\right)^{\frac 1 2} $  and $\|f \|_w = \left(\int_{\partial\Omega_w}| f|^2 ds\right)^{\frac 1 2} $ be the $L^2$ norm defined in the domain and on the domain boundary. The system \eqref{sys_nd}-\eqref{bd_nd} satisfies the following energy law.
\begin{theorem}\label{thm:edp}
	If $\phi$, $\vel$ and $P$ are  smooth solutions of the above system \eqref{sys_nd}-\eqref{bd_nd}, then  the following energy law is satisfied:
	\be\label{energye3}
	\frac{d}{dt}\mathcal{E}_{total}&\!\!\!=\!\!\!&\frac{d}{dt}(\mathcal{E}_{kin}+\mathcal{E}_{cell}+\mathcal{E}_w)\nonumber\\
	&=&  -2\|\eta^{1/2}\mathbf{D}_{\eta}\|^2  - \mathcal{M}\|\mu\|^2- \kappa \|\dot{\phi}\|_{w}^2-\|l_s^{-1/2} \vel_{\tau}\|^2_{w}~,
	\ee
	where 
	$\mathcal{E}_{total}=\mathcal{E}_{kin}+\mathcal{E}_{cell}+\mathcal{E}_{w}$,
	$\mathcal{E}_{kin}= \displaystyle\frac{1} 2 \int_{\Omega}|\vel|^2 d\bx $,
	$\mathcal{E}_{cell} =  \displaystyle\frac{\kappa_B }{2\epsilon} \int_{\Omega}|f|^2 d\bx+\mathcal{M}_v\frac{(V(\phi)-V(\phi_0))^2}{2V(\phi_0)}+\mathcal{M}_s\frac{(S(\phi)-S(\phi_0))^2}{2S(\phi_0)}$
	  and
	$\mathcal{E}_w\!\! =\!\! \alpha_w\!\!  \displaystyle\int_{\partial\Omega_w}f_w ds$. 
\end{theorem}

\noindent \textbf{Proof:} Multiplying the first equation in Eq.~\eqref{sys_nd}  with $\vel$ and integration by parts yield 
\be\label{kin_nd}
\frac{d}{dt}\mathcal{E}_{kin}&=&-\int_{\Omega}2\eta|\bD_{\eta}|^2d\bx+\int_{\partial\Omega_w}(\bsigma_{\eta}\cdot\bn)\cdot \vel_{\tau}ds +\int_{\Omega}\mu\nabla\phi\cdot\vel d\bx-\int_{\Omega}\lambda\delta_\epsilon \mathcal{P} : \nabla \vel d\bx\nonumber\\
&&+\int_{\partial\Omega_w}( \lambda\delta_\epsilon \mathcal{P}\cdot \bn) \cdot\vel_{\tau} ds\nonumber\\
&=&-\int_{\Omega}2\eta|\bD_{\eta}|^2d\bx-\int_{\Omega}\lambda\delta_\epsilon \mathcal{P} : \nabla \vel d\bx-l_s^{-1}\int_{\partial\Omega_w}|\vel_{\tau}|^2ds\nonumber\\ &&+\int_{\partial\Omega_w} L(\phi)\partial_{\tau}\phi\cdot\vel_{\tau}ds+\int_{\Omega}\mu\nabla\phi\cdot\vel d\bx~,
\ee 
where  the slip boundary condition in Eq.~\eqref{bd_nd} is applied.


Taking the inner product of the third equation in Eq.~\eqref{sys_nd}   with $\mu$ results in
 
\be\label{phit1}
\int_{\Omega}\frac{\partial\phi}{\partial t}\mu d\bx +\int_{\Omega}\vel\cdot\nabla\phi\mu d\bx =-\mathcal{M}\int_{\Omega}|\mu|^2d\bx~.
\ee
Multiplying the fourth equation in Eq.~\eqref{sys_nd}  with $\frac{\partial\phi}{\partial t}$ and integration by part give rise to 
\be\label{mut}
\int_{\Omega}\mu\frac{\partial\phi}{\partial t}d\bx &=&\kappa_B\int_{\Omega}g\frac{\partial\phi}{\partial t}d\bx +\frac{d}{dt}\left(\mathcal{M}_v\frac{(V(\phi)-V(\phi_0))^2}{2V(\phi_0)}\right)+\mathcal{M}_s\frac{S(\phi)-S(\phi_0)}{S(\phi_0)}\int_{\Omega}f\frac{\partial\phi}{\partial t}d\bx\nonumber\\
&=&\kappa_B\int_{\Omega}f\frac{\partial}{\partial t}\left(-\Delta \phi+\frac{1}{\epsilon^2}(\phi^3-\phi)\right)d\bx-\kappa_B\int_{\partial\Omega_w}\partial_nf\frac{\partial\phi}{\partial t}ds\nonumber\\
&&+\frac{d}{dt}\left(\mathcal{M}_v\frac{(V(\phi)-V(\phi_0))^2}{V(\phi_0)}\right)+\mathcal{M}_s\frac{d}{dt}\left(\frac{(S(\phi)-S(\phi_0))^2}{2S(\phi_0)}\right)\nonumber\\
&&-\mathcal{M}_s\left(\frac{S(\phi)-S(\phi_0)}{ S(\phi_0)}\right)\int_{\partial\Omega}\epsilon\partial_n\phi\frac{\partial\phi}{\partial t} ds\nonumber\\
&=&\frac{d}{dt}\left(\kappa_B\int_{\Omega}\frac{|f|^2}{2\epsilon}d\bx\right) +\frac{d}{dt}\left(\mathcal{M}_v\frac{(V(\phi)-V(\phi_0))^2}{2V(\phi_0)}\right)\nonumber\\
&&+\mathcal{M}_s\frac{d}{dt}\left(\frac{(S(\phi)-S(\phi_0))^2}{2S(\phi_0)}\right)-\int_{\partial\Omega} L(\phi) \frac{\partial\phi}{\partial t}ds +\alpha_w\frac{d}{dt}\int_{\partial\Omega_w}f_wds\nonumber\\
&=&\frac{d}{dt}(\mathcal{E}_{cell}+\mathcal{E}_w)-\int_{\partial\Omega_w} L(\phi) \frac{\partial\phi}{\partial t}ds~,
\ee
where  the definitions of $f(\phi)$, $g(\phi)$ and the boundary conditions of $\phi$ and $f$ are utilized. 

Multiplying the last equations with $\lambda$  and integration by parts leads to
\be\label{lamda}
&&  \int_\Omega (\lambda\delta_\epsilon \mathcal{P})  :\nabla\vel d\bx=0~.
\ee

Finally, the energy dissipation law \eqref{energye3} is obtained by combining Eqs. \eqref{kin_nd}, \eqref{phit1}, \eqref{mut} and \eqref{lamda}. 
$\null \hfill \blacksquare$

\section{Numerical Scheme and Discrete Energy law}

The numerical scheme for solving Eqs.~\eqref{sys_nd}-\eqref{bd_nd} uses the mid-point method for temporal discretization. Let $\Delta t$ denote the time step size, $()^{n+1}$ and $()^{n}$ denote the value of the variables at time $(n+1)\Delta t$ and $n\Delta t$. The semi-discrete equations in time are as follows:

\be
\left\{\begin{array}{ll}\label{sys_nd_dis}
	\frac{\vel^{n+1}-\vel^{n}}{\Delta t}+(\vel^{n+\frac{1}{2}}\cdot\na) \vel^{n+\frac{1}{2}}+\frac{1}{Re}\na P^{n+\frac{1}{2}}=\frac{1}{Re}\na\cdot(\eta^{n}(\na \vel^{n+\frac{1}{2}} +(\na \vel^{n+\frac{1}{2}})^T))&\\[3mm]
	+\frac{1}{Re}\mu^{n+\frac{1}{2}}\na\phi^{n+\frac{1}{2}}+\na\cdot \left(\lambda^{n+\frac{1}{2}} \mathcal{P}^{n} \delta_\epsilon\right)~,& \mbox{in~}\Omega~,\\[3mm]
	\na\cdot \vel^{n+\frac{1}{2}}=0~,& \mbox{in~}\Omega~,\\[3mm]
	\frac{\phi^{n+1}-\phi^{n}}{\Delta t}+(\vel^{n+\frac{1}{2}}\cdot\na)\phi^{n+\frac{1}{2}}=-\mathcal{M}\mu^{n+\frac{1}{2}}~,& \mbox{in~}\Omega~,\\[3mm]
	\mu^{n+\frac{1}{2}} =     \kappa_B g(\phi^{n+1},\phi^n)
	+\mathcal{M}_v\frac{(V(\phi^{n+\frac{1}{2}})-V(\phi_0))}{V(\phi_0)} 
	+\mathcal{M}_s\frac{(S(\phi^{n+\frac{1}{2}})-S(\phi_0))}{S(G_0)}f(\phi^{n+1},\phi^n)~,  & \mbox{in~}\Omega~,\\[3mm]
	f^{n+\frac{1}{2}} = -\epsilon \Delta\phi^{n+\frac{1}{2}}+\frac{1}{\epsilon}((\phi^{n+\frac{1}{2}})^2-1)\phi^{n+\frac{1}{2}}~, & \mbox{in~}\Omega~,\\[3mm]
\xi \epsilon^2 \nabla \cdot((\phi^{n})^2\nabla \lambda^{n+\frac1 2})+	\delta_\epsilon \mathcal{P}^{n} : \na\vel^{n+\frac{1}{2}}=0~, & \mbox{in~}\Omega~.
\end{array}\right.
\ee

The numerical boundary conditions can be written as:

\be\left\{\begin{array}{ll}\label{bd_nd_dis}
	\kappa\dot \phi^{n+\frac{1}{2}}=-L^{n+\frac{1}{2}}~,& \mbox{on~} \partial\Omega_w~,\\
	L^{n+\frac{1}{2}}= \kappa_B \partial_n f^{n+\frac{1}{2}}+ \mathcal{M}_s \epsilon\frac{S(\phi^{n+\frac{1}{2}})-S_0}{S_0}\partial_n\phi^{n+\frac{1}{2}}+ \alpha_w\frac{f_w^{n+1}-f_w^{n}}{\phi^{n+1}-\phi^{n}}~,  & \mbox{on~} \partial\Omega_w~,\\
	-l_s^{-1} u_{\tau_i}^{n+\frac{1}{2}} = \boldsymbol{\tau_i}\cdot (\eta^{n}(\na \vel^{n+\frac{1}{2}} +(\na \vel^{n+\frac{1}{2}})^T)+\lambda^{n+\frac{1}{2}}\delta_\epsilon \mathcal{P}^{n})\cdot \boldsymbol{n} \\
	- L^{n+\frac{1}{2}}\partial_{\tau_i}\phi^{n+\frac{1}{2}} ~,~i=1,2, & \mbox{on~}  \partial\Omega_w~,\\
  f^{n+\frac{1}{2}}=0~, &\mbox{on~} \partial\Omega_w~,\\
  \partial_n \lambda^{n+\frac 1 2} =0~,&\mbox{on~} \partial\Omega_w~,
\end{array}\right.\ee
where 
\be
f(\phi^{n+1},\phi^n) = -\epsilon \Delta \phi^{n+\frac{1}{2}}+\frac{1}{4\epsilon}((\phi^{n+1})^2+(\phi^{n})^2-2)(\phi^{n+1}+\phi^n)~,\\
g(\phi^{n+1},\phi^n) = \left(-\Delta f^{n+\frac{1}{2}}+\frac{1}{\epsilon^2}\left((\phi^{n+1})^2+(\phi^{n})^2+\phi^{n+1}\phi^{n}-1\right)f^{n+\frac{1}{2}}\right)~, 
\ee
$(\cdot)^{n+\frac 1 2} = \frac {(\cdot)^n +(\cdot)^{n+1}}{2}$ and $\mathcal {P}^n = I -\bn^{n}_m\otimes\bn^{n}_m$ with $\bn_m^n = \frac{\nabla\phi^{n}}{|\nabla\phi^n|}$.
\begin{rmk}
Here  an extra term $\xi \epsilon^2 \nabla \cdot((\phi^{n})^2\nabla \lambda^{n+\frac1 2})$ is introduced  to extend $\lambda$ to the whole domain as in \cite{voigt_2014_local_inextensibility}. This can be regarded as a relaxation for  the local inextensibility near the membrane. The last equation can also be derived directly by the variational method described in Section 2 by adding an extra dissipation term $\int_{\Omega} \xi \epsilon^2 |\phi\nabla\lambda|^2d\bx$ in the dissipation functional. 
\end{rmk}

\begin{rmk}
 $\delta_\epsilon$ is approximated by $\frac{1}{\sigma}|\nabla \phi^{n}|^2$, where $\sigma$ is a positive constant that adjusts the value of $\delta_\epsilon$ to be close to 1 at $\phi =0$, center of the diffuse interface \cite{lee2012regularized}.
\end{rmk}

The above scheme obeys the following  theorem of energy stability. 
\begin{theorem}
\label{energyTh_dis}
    If $(\phi^n,\vel^n,P^n)$ are  smooth solutions of the above system \eqref{sys_nd_dis}-\eqref{bd_nd_dis}, then  the following energy law is satisfied:
	\be\label{energye3_dis}
 \mathcal{E}_{total}^{n+1}-\mathcal{E}_{total}^n&\!\!\!=\!\!\!& (\mathcal{E}_{kin}^{n+1}+\mathcal{E}_{cell}^{n+1}+\mathcal{E}_w^{n+1})- (\mathcal{E}_{kin}^{n}+\mathcal{E}_{cell}^{n}+\mathcal{E}_w^{n})\nonumber\\
	&=&  \triangle t\left(-2\|(\eta^n)^{1/2}\mathbf{D}^{n+\frac 12}_{\eta}\|^2  - \mathcal{M}\|\mu^{n+\frac 1 2}\|^2- \xi\|\ \epsilon \phi^n \nabla \lambda^{n+\frac 1 2}  \|^2\right.\nonumber\\
	&&\left.- \frac 1 {\kappa} \|L(\phi^{n+\frac 1 2})\|_{w}^2-\|l_s^{-1/2} \vel_{\tau}^{n+\frac 1 2}\|^2_{w}\right)~,
	\ee
	where 
	$\mathcal{E}_{total}^n=\mathcal{E}_{kin}^n+\mathcal{E}_{cell}^n+\mathcal{E}_{w}^n$ with
	$\mathcal{E}_{kin}^n=\frac{Re} 2 \|\vel^n\|^2$,
	$\mathcal{E}_{cell}^n = \frac{\kappa_B\|f^n\|^2}{2\epsilon}+\mathcal{M}_v\frac{(V(\phi^n)-V(\phi_0))^2}{2V(\phi_0)}+\mathcal{M}_s\frac{(S(\phi^ n)-S(\phi_0))^2}{2S(\phi_0)}$
	  and
	$\mathcal{E}_w^n\!\! =\!\! \alpha_w\!\!  \displaystyle\int_{\partial\Omega_w}f_w^n ds$.

\end{theorem}

The following two important lemmas are needed for proving Theorem \ref{energyTh_dis}.
\begin{lemma}\label{lma:1}
Let
\be 
f(\phi^{n+1},\phi^{n})=-\epsilon \Delta \phi^{n+\frac{1}{2}}+\frac{1}{4\epsilon}((\phi^{n+1})^2+(\phi^{n})^2-2)(\phi^{n+1}+\phi^n)~.
\ee
Then $f(\phi^{n+1},\phi^{n})$ satisfies
\be
\int_\Omega f(\phi^{n+1},\phi^{n})(\phi^{n+1}-\phi^{n})d\bx=S^{n+1}-S^n-\int_{\partial\Omega_w} \epsilon\partial_n \phi^{n+\frac{1}{2}}(\phi^{n+1}-\phi^{n})ds~,
\ee
where $S^{n+1}=\int_\Omega G(\phi^{n+1})d\bx, S^{n}=\int_\Omega G(\phi^{n})d\bx$.

\noindent\textbf{Proof:}
\be
&&\int_\Omega f(\phi^{n+1},\phi^{n})(\phi^{n+1}-\phi^{n})d\bx\nonumber\\
&=&\int_\Omega\left( \epsilon\Delta\phi^{n+\frac{1}{2}}(\phi^{n+1}-\phi^{n})+\frac{1}{4\epsilon}((\phi^{n+1})^2+(\phi^{n})^2-2)(\phi^{n+1}+\phi^n)(\phi^{n+1}-\phi^{n})\right)d\bx \nonumber\\
&=&\int_\Omega \epsilon\na\phi^{n+\frac{1}{2}}\cdot\na (\phi^{n+1}-\phi^{n}) d\bx -\int_{\partial\Omega_W}\epsilon\partial_n\phi^{n+\frac{1}{2}}(\phi^{n+1}-\phi^{n})ds\nonumber\\
&&+\int_\Omega \frac{1}{4\epsilon}((\phi^{n+1})^4-2(\phi^{n+1})^2-(\phi^{n})^4+2(\phi^{n})^2)d\bx\nonumber\\
&=&\int_\Omega \left[\frac{\epsilon}{2}((\na\phi^{n+1})^2-(\na\phi^n)^2)+\frac{1}{4\epsilon}(((\phi^{n+1})^2-1)^2-((\phi^{n})^2-1)^2)\right]d\bx\nonumber\\
&&-\int_{\partial\Omega_W}\epsilon\partial_n\phi^{n+\frac{1}{2}}(\phi^{n+1}-\phi^{n})ds\nonumber\\
&=&\int_\Omega\left( \frac{\epsilon}{2}(\na\phi^{n+1})^2+\frac{1}{4\epsilon}(\phi^{n+1}-1)^2\right)d\bx-\int_\Omega\left( \frac{\epsilon}{2}(\na\phi^{n})^2+\frac{1}{4\epsilon}(\phi^{n}-1)^2\right)d\bx\nonumber\\
&&-\int_{\partial\Omega_w}\epsilon\partial_n\phi^{n+\frac{1}{2}}(\phi^{n+1}-\phi^{n})ds\nonumber\\
&=&S^{n+1}-S^n-\int_{\partial\Omega_w} \epsilon\partial_n \phi^{n+\frac{1}{2}}(\phi^{n+1}-\phi^{n})ds~. \nonumber
\ee
$\hfill\blacksquare$
\end{lemma}

\begin{lemma}\label{lma:2}
Let
\be
g(\phi^{n+1},\phi^{n})=-\Delta f^{n+\frac{1}{2}}+\frac{1}{\epsilon^2}((\phi^{n+1})^2+(\phi^{n})^2+\phi^{n+1}\phi^n-1)f^{n+\frac{1}{2}}~.
\ee
Then $g(\phi^{n+1},\phi^{n})$ satisfies
\be
\int_\Omega g(\phi^{n+1},\phi^n)(\phi^{n+1}-\phi^n)d\bx=\int_\Omega \frac{1}{2\epsilon}((f^{n+1})^2-(f^n)^2)d\bx-\int_{\partial\Omega}\partial_n f^{n+\frac{1}{2}}(\phi^{n+1}-\phi^n)ds~,
\ee
where $f^{n+1}=-\epsilon \Delta \phi^{n+1}+\frac{1}{\epsilon}((\phi^{n+1})^2-1)\phi^{n+1}, f^{n}=-\epsilon \Delta \phi^{n}+\frac{1}{\epsilon}((\phi^{n})^2-1)\phi^{n}$.

\noindent\textbf{Proof:}
\be
&&\int_\Omega g(\phi^{n+1},\phi^n)(\phi^{n+1}-\phi^n)d\bx \nonumber\\
&=&\int_\Omega \na f^{n+\frac{1}{2}}\na (\phi^{n+1}-\phi^n)d\bx-\int_{\partial\Omega} \partial_n f^{n+\frac{1}{2}}(\phi^{n+1}-\phi^n)ds\nonumber\\
&&+\int_\Omega \frac{1}{\epsilon^2}f^{n+\frac{1}{2}}(((\phi^{n+1})^2-1)\phi^{n+1}-((\phi^{n})^2-1)\phi^n)d\bx\nonumber\\
&=&-\int_\Omega f^{n+\frac{1}{2}}\Delta (\phi^{n+1}-\phi^n)d\bx-\int_{\partial\Omega_w} \partial_n f^{n+\frac{1}{2}}(\phi^{n+1}-\phi^n)ds \nonumber\\
&&+\int_\Omega \frac{1}{\epsilon^2}f^{n+\frac{1}{2}}(((\phi^{n+1})^2-1)\phi^{n+1}-((\phi^{n})^2-1)\phi^n)d\bx\nonumber\\
&=&\int_\Omega \frac{1}{\epsilon}f^{n+\frac{1}{2}}\left((-\epsilon\Delta\phi^{n+1}+\frac{1}{\epsilon}((\phi^{n+1})^2-1)\phi^{n+1})-(-\epsilon\Delta\phi^{n}+\frac{1}{\epsilon}((\phi^{n})^2-1)\phi^{n})\right)d\bx\nonumber\\
&&-\int_{\partial\Omega_w} \partial_n f^{n+\frac{1}{2}}(\phi^{n+1}-\phi^n)ds\nonumber\\
&=&\int_\Omega \frac{1}{2\epsilon}(f^{n+1}+f^{n})(f^{n+1}-f^{n})d\bx-\int_{\partial\Omega_w} \partial_n f^{n+\frac{1}{2}}(\phi^{n+1}-\phi^n)ds\nonumber\\
&=&\int_\Omega \frac{1}{2\epsilon}((f^{n+1})^2-(f^{n})^2)d\bx-\int_{\partial\Omega_w} \partial_n f^{n+\frac{1}{2}}(\phi^{n+1}-\phi^n)ds\nonumber~.
\ee
$\hfill\blacksquare$
\end{lemma}

\noindent{\textbf{Proof of Theorem \ref{energyTh_dis}}:}  Multiplying the first equation in system \eqref{sys_nd_dis} by $\Delta t \vel^{n+\frac{1}{2}}$ gives
\be\label{kenetic_dis}
&&\int_{\Omega} \frac{1}{2}((\vel^{n+1})^2-(\vel^{n})^2)d\bx+\int_{\Omega}\Delta t \vel^{n+\frac{1}{2}}\cdot ((\vel^{n+\frac{1}{2}}\na)\cdot \vel^{n+\frac{1}{2}})d\bx-\frac{\Delta t}{Re}\int_{\Omega}P^{n+\frac{1}{2}}\na\cdot \vel^{n+\frac{1}{2}}d\bx\nonumber\\
&=&-\frac{\Delta t}{Re}\int_{\Omega}\na \vel^{n+\frac{1}{2}}:\eta^{n}(\na \vel^{n+\frac{1}{2}}+(\na \vel^{n+\frac{1}{2}})^T)d\bx+\frac{\Delta t}{Re}\int_{\Omega}\vel^{n+\frac{1}{2}}\cdot\na \phi^{n+1}\mu ^{n+1}d\bx\nonumber\\
&&-\Delta t\int_{\Omega}\lambda\delta_\epsilon \mathcal{P}^n : \nabla \vel^{n+\frac 1 2} d\bx+\frac{\Delta t}{Re}\int_{\partial\Omega_w} \lambda^{n+\frac{1}{2}}( \delta_\epsilon\mathcal{P}^{n}\cdot\bn)\cdot \vel_\tau ^{n+\frac{1}{2}} ds\nonumber\\
&&+\frac{\Delta t}{Re}\int_{\partial\Omega_w}\vel^{n+\frac{1}{2}}\cdot\eta^n((\na \vel^{n+\frac{1}{2}}+(\na\vel^{n+\frac{1}{2}})^T)\cdot \boldsymbol{n})ds~.
\ee

Multiplying the fourth equation in system \eqref{sys_nd_dis}  by $\frac{\phi^{n+1}-\phi^{n}}{Re}$ and integration by parts lead to

\be\label{mu_dis}
\frac{1}{Re}\int_\Omega \mu^{n+1}(\phi^{n+1}-\phi^{n})d\bx=\frac{\kappa_B}{Re}\int_\Omega \frac{1}{2\epsilon}((f^{n+1})^2-(f^{n})^2)d\bx\nonumber\\
+\frac{\mathcal{M}_v}{Re}\frac{(V(\phi^{n+1})-V_0)^{2}-(V(\phi^{n})-V_0)^{2}}{2V_0}
+\frac{\mathcal{M}_s}{Re}\frac{(S(\phi^{n+1})-S_0)^{2}-(S(\phi^{n})-S_0)^{2}}{2S_0}\nonumber\\
-\int_{\partial\Omega} \partial_n f^{n+\frac{1}{2}}(\phi^{n+1}-\phi^{n})ds
-\frac{\mathcal{M}_s}{Re}\int_{\partial\Omega_w} \frac{S(\phi^{n+\frac{1}{2}})-S_0}{S_0}\epsilon\partial_n\phi^{n+\frac{1}{2}}(\phi^{n+1}-\phi^{n})ds~.
\ee

Multiplying the third equation in system \eqref{sys_nd_dis}  by $\frac{\mu^{n+1}\Delta t}{Re}$ and integration by parts yield

\be\label{phit_dis}
\frac{1}{Re}\int_\Omega \mu^{n+1}(\phi^{n+1}-\phi^{n})d\bx+\frac{\Delta t}{Re}\int_{\Omega}\mu^{n+1}(\vel^{n+1/2}\cdot\na)\phi^{n+1}d\bx=-\frac{\mathcal{M}\Delta t}{Re}\int_\Omega (\mu^{n+1})^2d\bx~.
\ee
Multiplying the last equation in system \eqref{sys_nd_dis}  by $\lambda^{n+\frac 1 2}\Delta t$ and integration by parts give 

\be\label{lamda_dis}
&&-\Delta t\int_\Omega \xi \epsilon^2 (\phi^n)^2 \left|\nabla \lambda^{n+\frac 1 2}  \right|^2 d\bx + \Delta t\int_\Omega (\lambda^{n+\frac{1} 2}\delta_\epsilon \mathcal{P}^n)  :\nabla\vel^{n+\frac 1 2} d\bx=0~.
\ee
The discretized energy dissipation law \eqref{energye3_dis} is obtained by combining Eqs.~\eqref{kenetic_dis}-\eqref{lamda_dis}.
$\hfill \blacksquare$

\section{Simulation Results}
Numerical simulations using the model introduced in the paper are presented in this section. The first example is used to illustrate the convergence and energy stability of the proposed numerical scheme. Then feasibility of the proposed model and the model simulation scheme  to studying vesicle motion and shape transformation is assessed by cell tank treading and tumbling tests. The last simulation is devoted to studying   effects  of mechanical and geometric properties of a vesicle on its deformability when it passes through a narrow channel.  

\subsection{Convergence study}
\begin{figure}[ht]
		\centering
		\includegraphics[width=4.in]{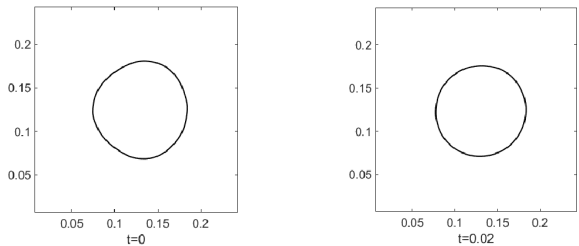}
		\caption{Relaxation of a tear shape vesicle.}
		\label{fig:convergence_space}
	\end{figure}

The initial condition of the convergence test is set to be a 2D tear shape vesicle in a closed cube with  intercellular and extracellular fluid velocity being 0. Thanks to the bending force of the cell membrane, the shape of the vesicle  gradually transforms into a perfect circle to minimize the total energy (see Figure \ref{fig:convergence_space}).
The parameter values used for  this simulation are chosen as follows:
$Re=2\times 10^{-4}$,   $\mathcal{M}=5\times 10^{-5}$, $\kappa_B=8\times 10^{-1}$, $\epsilon =2.5\times 10^{-2}$, $\mathcal{M}_v=20$, $\mathcal{M}_s=2$, $\xi=1.6\times10^{5}$, $\kappa=8\times 10^{-10}$, $\alpha_w=2\times 10^{9}$, $l_s=5\times 10^{-3}$.

In the simulations, the numerical solution computed with a mesh size  $h=1/240$  is treated as the reference solution or ``the true solution''.  As shown in Table \ref{tab:convergence_space}, our scheme is a   second-order accurate in space. 
	\begin{table}[ht]
		\centering
		\begin{tabular}{|c|c|c|c|c|c|c|}
			\hline
			\makecell{Spacial mesh \\ size
				$h$} & \multicolumn{6}{|c|}{P2 Element} \\
			\cline{2-7}
			& Err($u_x$) & \makecell{ Convergence \\ Rate($u_x$)} &  Err($u_y$) & \makecell{ Convergence\\ Rate($u_y$)}& Err($\phi$) & \makecell{Convergence\\ Rate($\phi$)} \\ 
			\hline
			1/47& 1.3e-1 & & 1.5e-1 & & 1.4e-2 & \\ 
			\hline
			1/71& 8.3e-2 & 1.15& 7.6e-2 & 1.71& 6.1e-3 &1.97 \\ 
			\hline
			1/107& 3.8e-2 & 1.94 & 3.7e-2 & 1.83 & 2.3e-3 & 2.45 \\ 
			\hline
			1/160& 1.5e-2 & 2.35  & 1.3e-2 & 2.59 & 5.7e-4 & 3.42 \\ 
			\hline
		
		\end{tabular}
		\caption{$L^2$ norm of the error and convergence rate for  velocity $\vel=(u_x,u_y)$, phase-field function $\phi$,  at time $t=0.02 $ with both intercellular and extracellular fluid viscosity  being 1. }
		\label{tab:convergence_space}
	\end{table}
	
The time convergence rate of the scheme is obtained by comparing the numerical errors calculated using  each pair of successively reduced time step sizes. The purpose of doing so is to eliminate the influence from the error of the reference solution which is also a numerical result. Larger Reynolds number $Re$ and interface thickness $\epsilon$, and a smoother initial profile of the interface are applied to ensure that the convergence rate is not affected by any sharp changes in the phase field label function $\phi(\bx)$. Results in Table \ref{tab:convergence_time} confirm that our scheme is a also second-order accurate in time.  

\begin{figure}[!ht]
		\centering
		\includegraphics[width=5.in]{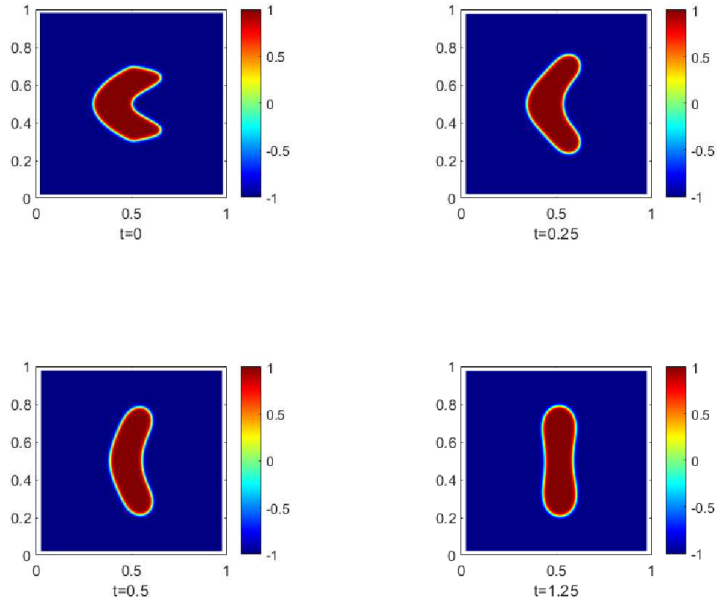}
		\caption{Relaxation of a bent vesicle. The fluid viscosities  are 1 and 50 for intercellular and extracellular fluids, respecfively.}
		\label{fig:profile}
	\end{figure}

	\begin{table}[ht]
		\centering
		\begin{tabular}{|c|c|c|c|c|c|c|}
			\hline
			{time step
				$\Delta t$} & \multicolumn{6}{|c|}{P2 Element} \\
			\cline{2-7}
			& Err($u_x$) & \makecell{ Convergence \\ Rate($u_x$)}&  Err($u_y$) & \makecell{ Convergence \\ Rate($u_y$)}& Err($\phi$) &  \makecell{ Convergence \\ Rate($\phi$)} \\ 
		    \hline
			0.025& - & & - & & -  & \\ 
			\hline
			0.0125& 8.12e-6 & & 8.13e-6 & &9.92e-6 & \\ 
			\hline
			0.00625& 2.90e-6 & 1.49& 2.97e-6 & 1.45& 2.42e-6 &2.04 \\ 
			\hline
			0.003125&1.03e-6  & 1.48 &1.07e-6  & 1.48 &5.98e-7  & 2.01 \\ 
			\hline
			0.0015625&2.53e-7  & 2.03 &2.60e-7  & 2.03 &1.49e-7  & 2.01 \\ 
			\hline
		
		\end{tabular}
		\caption{$L^2$ norm of the error and convergence rate for  velocity $\vel=(u_x,u_y)$, phase-field function $\phi$,  at time $t=0.05 $ with both intercellular and extracellular fluid viscosities  being 1. }
		\label{tab:convergence_time}
	\end{table}
Finally, the energy law (\textbf{ Theorem \ref{energyTh_dis}}) and conservation of  mass and surface area  of vesicles are tested by the relaxation of a bent vesicle, which gradually evolves back to its equilibrium biconcave shape. Figure \ref{fig:profile} shows the snapshots of the vesicle profile at different times $t= 0, 0.25, 0.5$ and $1.25$.  The parameter values used here are:

$Re=2\times 10^{-4}$,   $\mathcal{M}=2.5\times 10^{-3}$, $\kappa_B=2$, $\epsilon =7.5\times 10^{-3}$, $\mathcal{M}_v=20$, $\mathcal{M}_s=2$, $\xi=7.1\times 10^{4}$, $\kappa=2\times 10^{-10}$, $\alpha_w=2\times 10^{9}$, $l_s=0.5$.

\begin{figure}[!ht]
     \centering
     \begin{subfigure} 
         \centering
         \includegraphics[width=2.in]{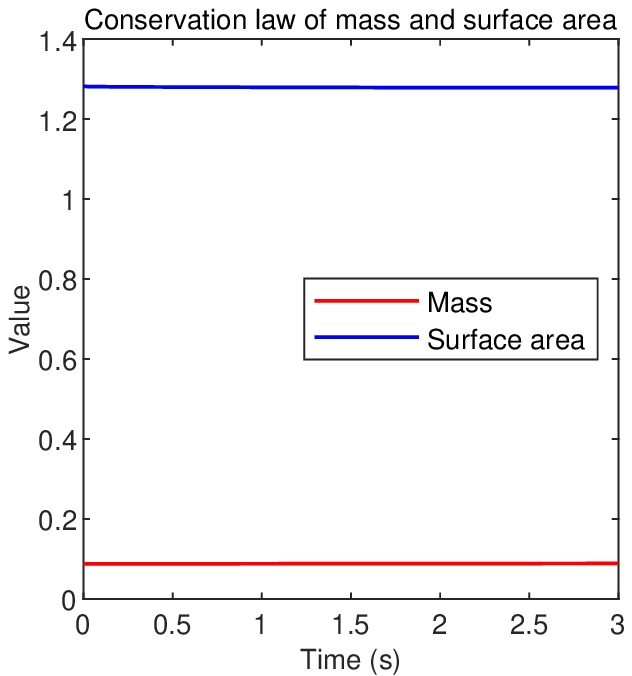}
     \end{subfigure}
     \begin{subfigure} 
         \centering
         \includegraphics[width=2.in]{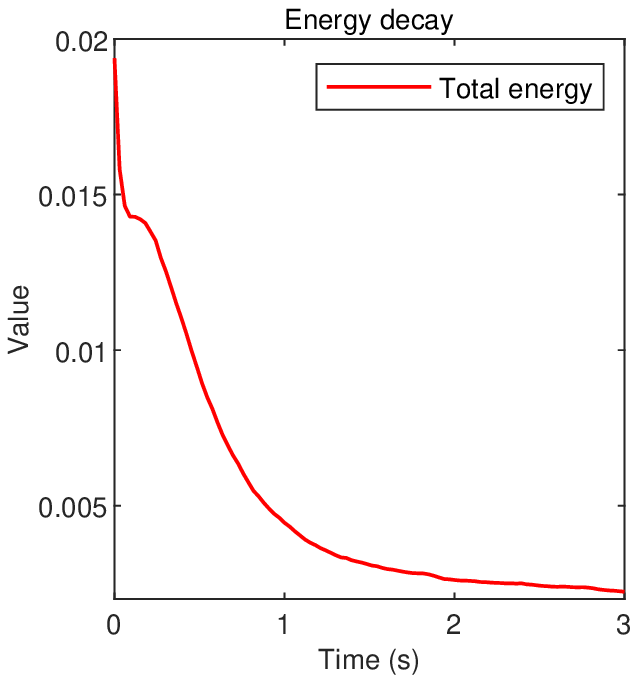}
     \end{subfigure}
        \caption{The test case of relaxation of a bent vesicle. Left: Change of mass and surface area vs. time; Right: Change of discrete energy vs. time.}
        \label{fig:conservation}
\end{figure}

The changes of vesicle mass and surface area and the change of total discrete energy of this test case computed by the scheme (Eqs.~\eqref{sys_nd}-\eqref{bd_nd}) are shown in Figure \ref{fig:conservation}. It is evident that the vesicle mass and surface area are almost perfectly preserved, and  the total energy decays over the course of time as expected.

\subsection{Tank treading and tumbling}
The vesicle motion in a Couette flow changes with respect to the ratio of the viscosities $\eta_{in}$ and  $\eta_{out}$ of  intracellular and extracellular  fluids \cite{laadhari2012vesicle,basu2011tank,fischer2004shape,hao2015fictitious}.  When this viscosity ratio is small, the vesicle is prone to move in the tank treading mode; while the tumbling mode is preferred when the viscosity ratio is large.   The parameter values utilized for this vesicle motion simulation are set as follows:

$Re=2\times 10^{-4}$, $\delta_{\epsilon}=|\nabla\phi^{n}|^{2}$, $\mathcal{M}=10^{-3}$, $\kappa_B=5\times 10^{-3}$, $\epsilon =7.5\times 10^{-3}$, $\mathcal{M}_v=20$, $\mathcal{M}_s=200$, $\xi=1.78\times 10^{7}$, $\kappa=2\times 10^{-12}$, $\alpha_w=2\times 10^{9}$, $l_s=0.2$.

\begin{figure}[h]
		\centering
		\includegraphics[width=4.in]{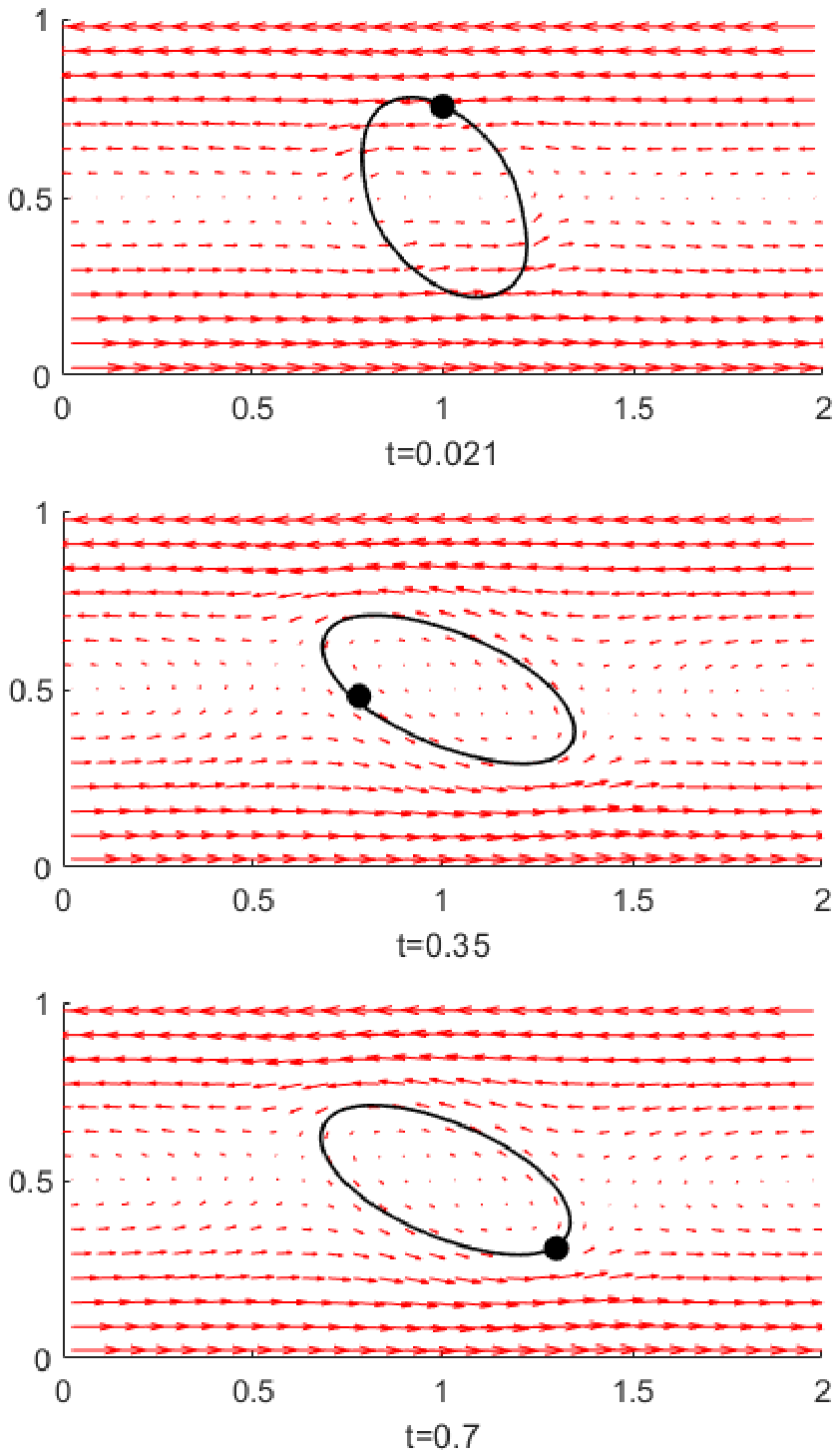}
		\caption{Tank treading with viscosity ratio $1:1$. The orientation of the vesicle and the velocity field are kept stable when the system comes to equilibrium. The tracer point (in black) on the membrane rotates in the clockwise direction along  the vesicle.}
		\label{fig:tank treading}
	\end{figure}

\begin{figure}[h]
		\centering
		\includegraphics[width=4.in]{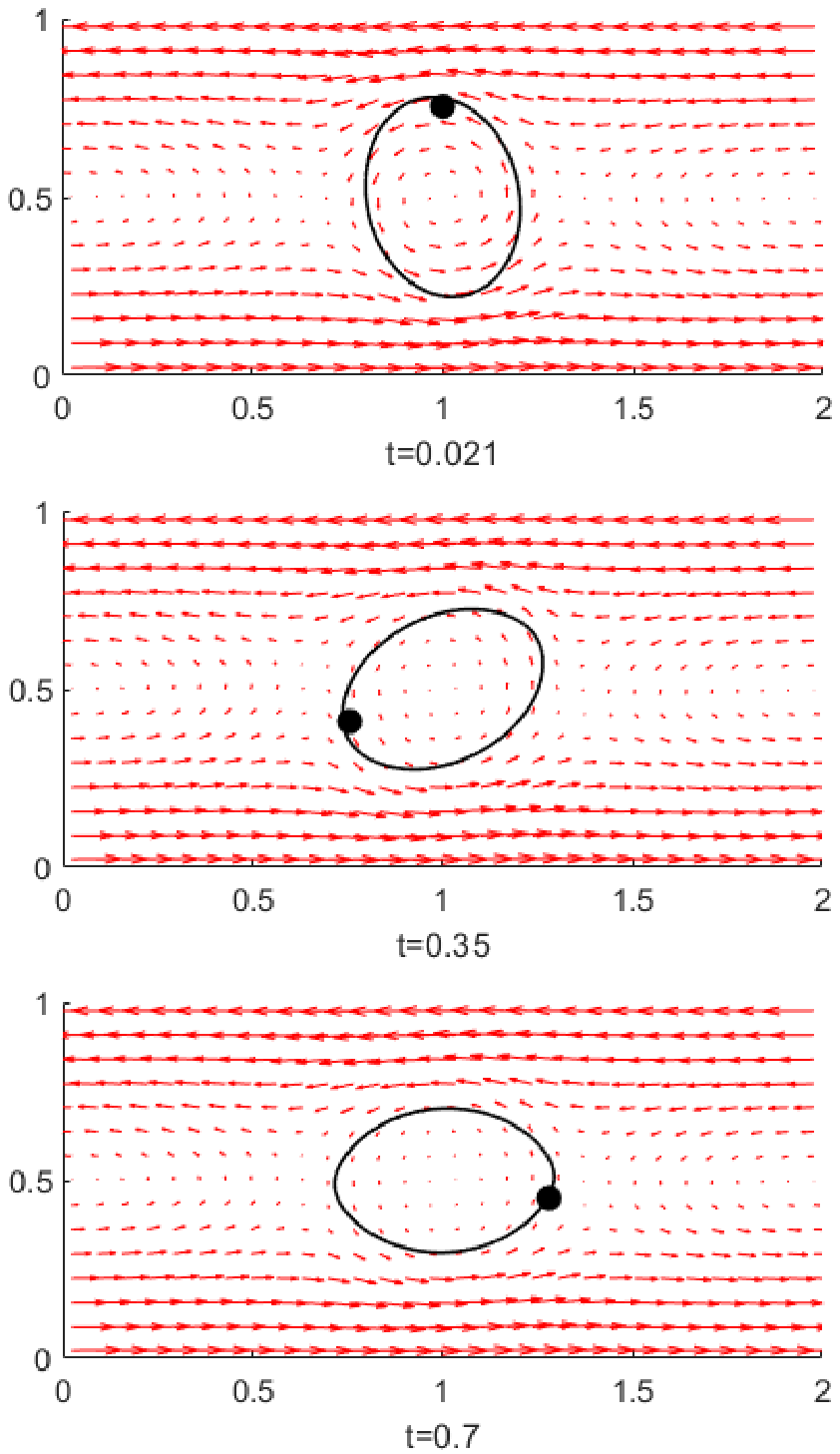}
		\caption{Tumbling with viscosity ratio 1:500. The vesicle keeps rotating in the flow. Position of the tracer point (in black) is fixed with respect to the vesicle membrane.}
		\label{fig:tumbling}
	\end{figure}

The  upper and bottom walls of the domain are set to move in opposite direction horizontally with velocities $-20$ and $20$, respectively. The simulation domain is  $2\times 1$, and the initial shape of the vesicle is chosen to be an ellipse with eccentricity $\sqrt{3}$.
The ratios of viscosities of the  intracellular and extracellular fluids are set to be $1:1$ and $1:500$, respectively. 
Figures \ref{fig:tank treading} and \ref{fig:tumbling} show the interfaces   of tank treading vesicle (low viscosity ratio case) and tumbling vesicle (high viscosity ratio case) and corresponding fluid velocity fields at different times, respectively.
A  point on the interface (black solid) is  tracked to illustrate  these two different types of motion. 
 For the tank treading motion, the angle between the long axis of the vesicle and horizontal axis is fixed when the vesicle is at equilibrium, but the tracer point rotates in a    counter clockwise  direction along the membrane. For the tumbling motion, the  vesicle keeps rotating and the tracer point does not move with respect to the membrane shape.    

\begin{figure}[h]
		\centering
		\includegraphics[width=3.in]{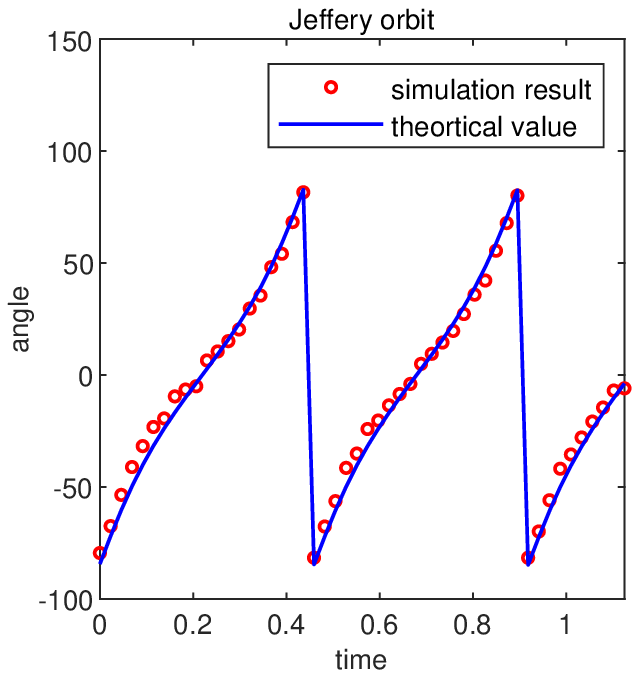}
		\caption{Comparison between theoretical and simulation results of the flipping ellipse. The blue line is the angle between the long axis of the ellipse and the horizontal axis predicted by the Jeffery orbit theory, and the red circles are the angle from the simulation.}
		\label{fig:Jeffery}
	\end{figure}

Next, simulation result of tumbling motion of a rigid ellipse is compared with the theoretical solution  obtained using Jeffery’s orbit theory \cite{jeffery1922motion}.  Specifically, the angle between the long axis of the ellipse and the horizontal axis is compared.  As shown in Figure \ref{fig:Jeffery}, our simulation result is in close agreement with the analytical Jeffery orbit.  

\subsection{Vesicle passing through a narrow fluid channel}


Finally,  the calibrated model is used to study  the effects of  mechanical properties of the membrane of the vesicle on its circulating through constricting micro channels \cite{han2019flow}. The vesicle shape is described by an ellipse with eccentricity $\sqrt{3}$, and   the width of the squeezing section of the narrow channel is 0.3 by default.  A pressure drop boundary condition is applied at the inlet(left) and outlet(right) of the domain by setting the pressure on the inlet and outlet to be $P=50$ and $P=-50$, respectively. Fluid viscosity ratio is set to be $1 :10$ for extracellular and intracellular fluids, respectively.
The other parameters  are  as follows:

$Re=2\times 10^{-4}$, $\delta_{\epsilon}=10\times|\nabla\phi^{n}|^{2}$, $\mathcal{M}=5\times10^{-4}$, $\kappa_B=4\times 10^{-2}$, $\epsilon =7.5 \times 10^{-3}$, $\mathcal{M}_v=20$, $\mathcal{M}_s=100$, $\xi=7.1\times 10^{4}$, $\kappa=4\times 10^{-11}$, $\alpha_w=2\times 10^{9}$, $l_s=5\times 10^{-3}$.

\begin{figure}[h]
\centering     
\subfigure [] {\includegraphics[width=2.5in]{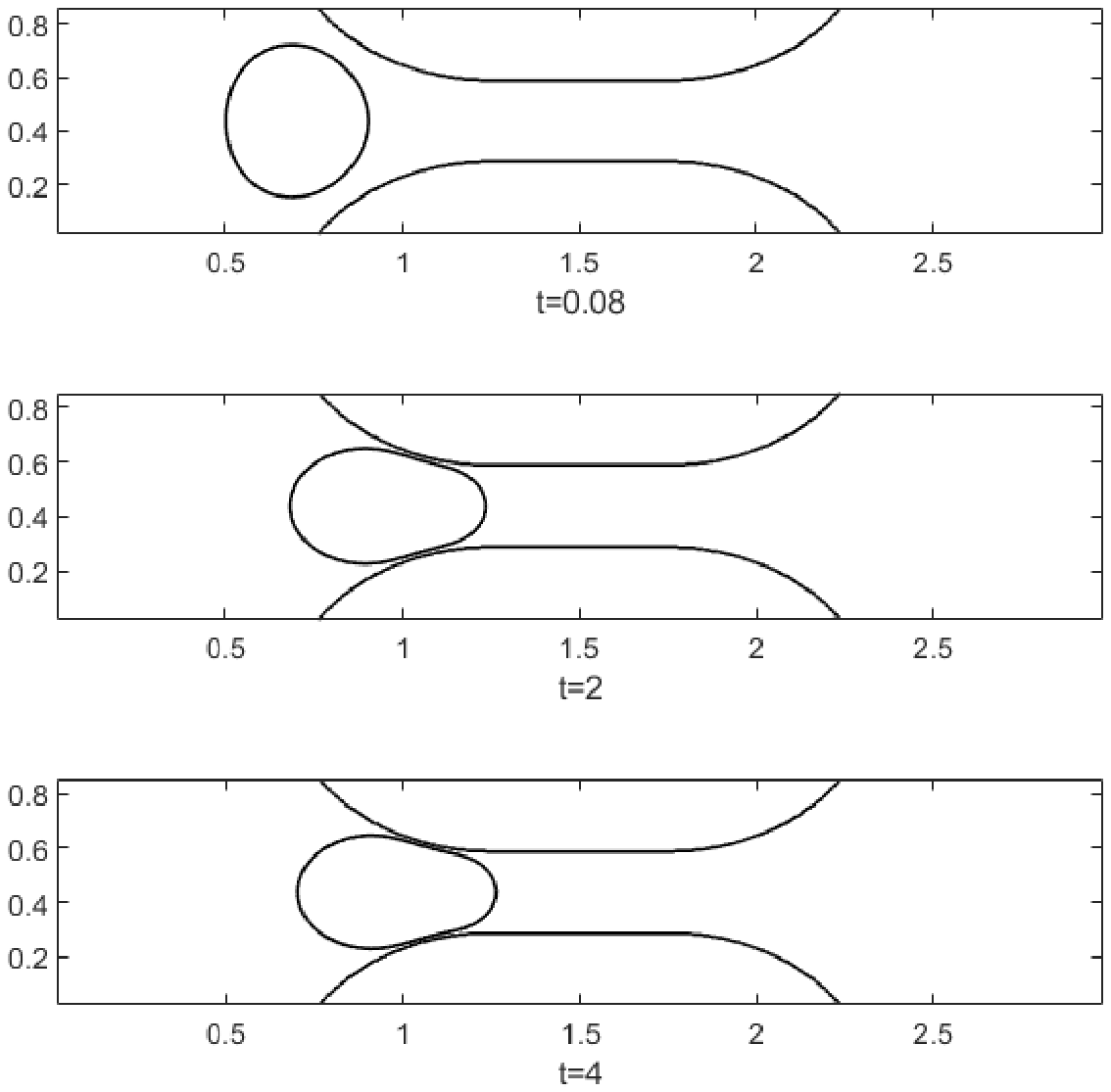}}
\subfigure[] { \includegraphics[width=2.5in]{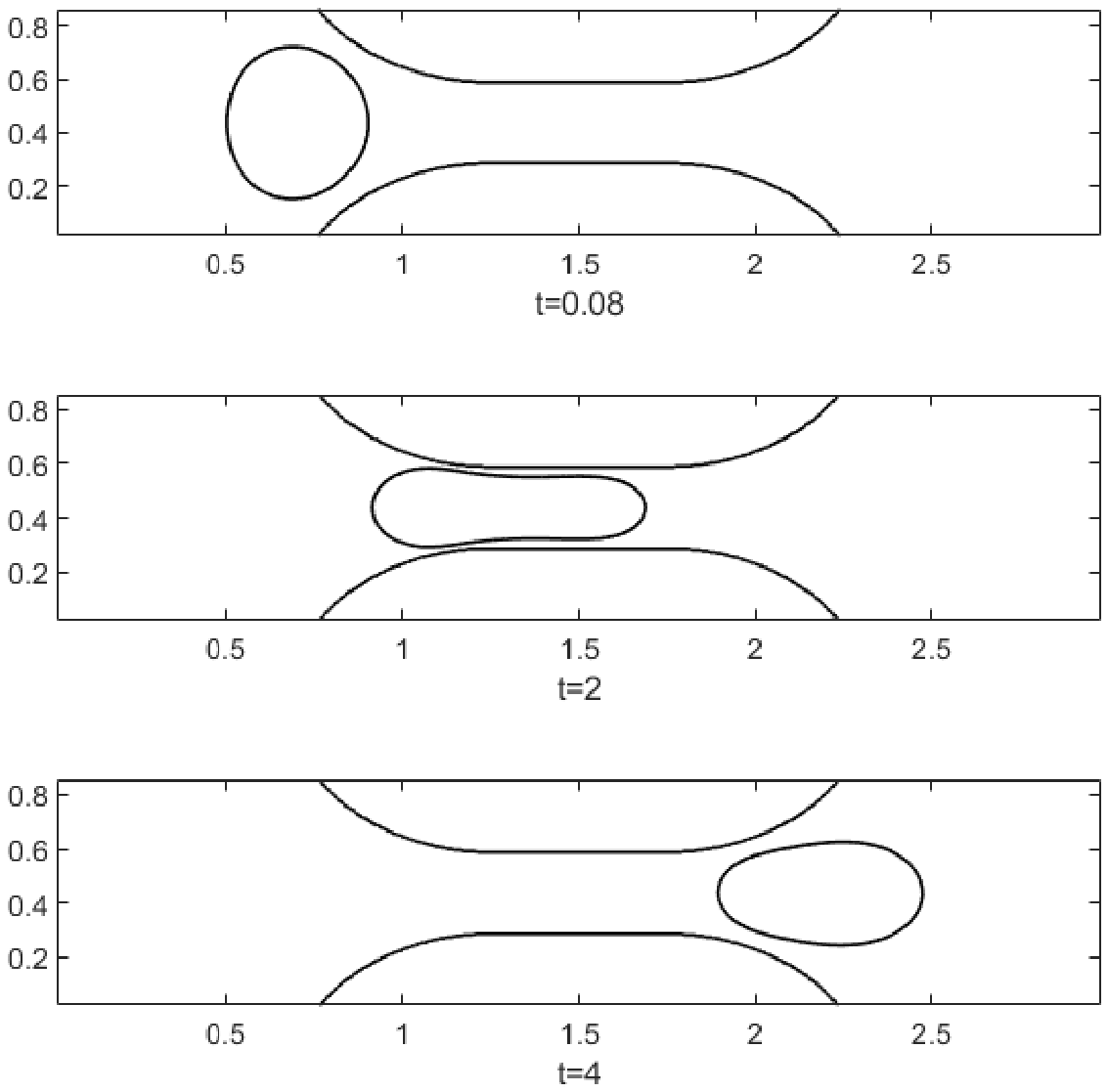}}
\subfigure[] { \includegraphics[width=2.5in]{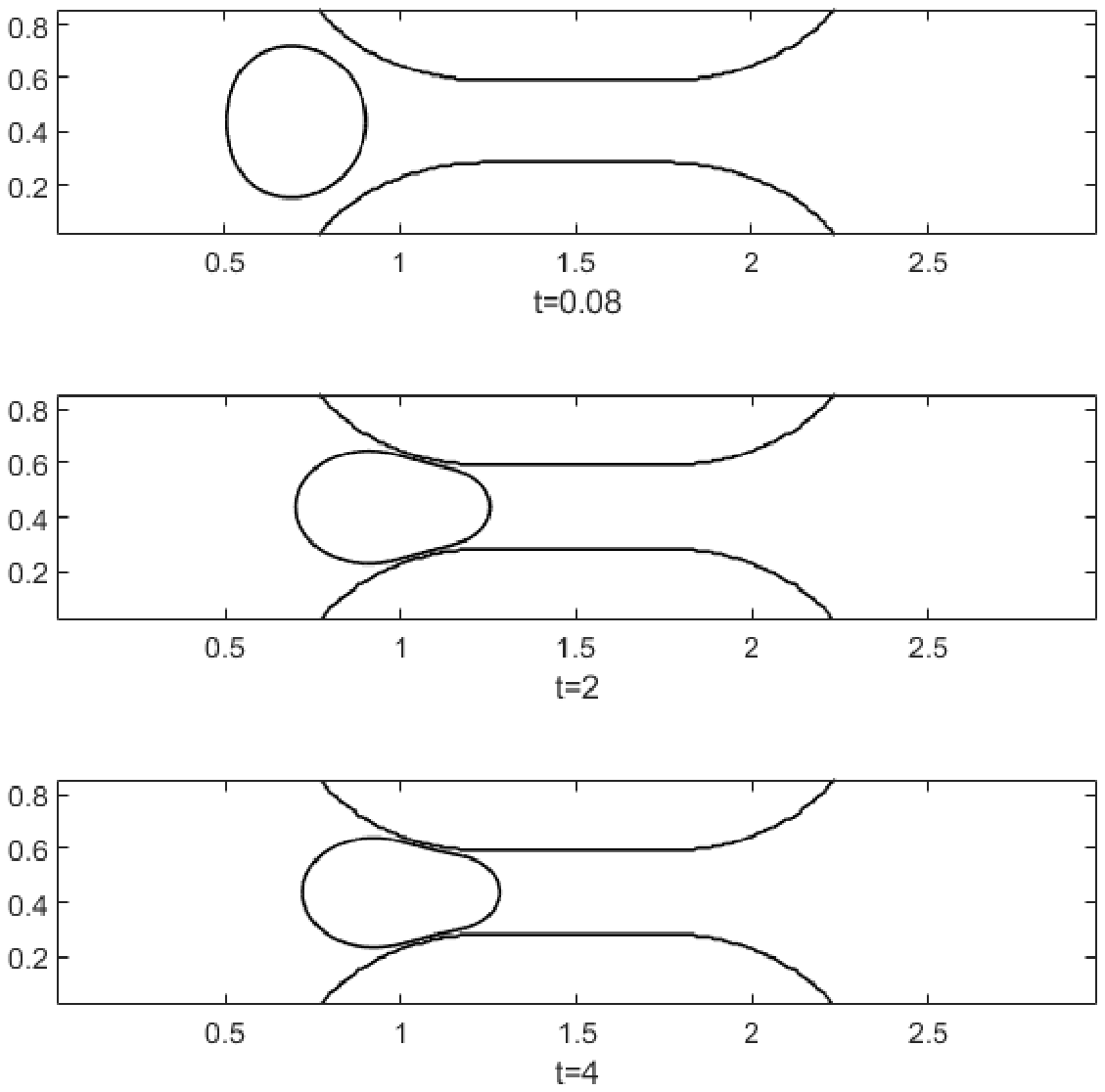}}
\caption{\label{fig:contour} Snapshots of vesicles with different surface area constraints at times $t = 0.08,2$ and $4$, respectively.  (a) $\mathcal{M}_s=100$ with the local inextensibility; (b) $\mathcal{M}_s=100$  without the local inextensibility;  (c) $\mathcal{M}_s=20000$ without the local inextensibility. }
\end{figure}

\begin{figure}[h]
		\centering
		\includegraphics[width=3.in]{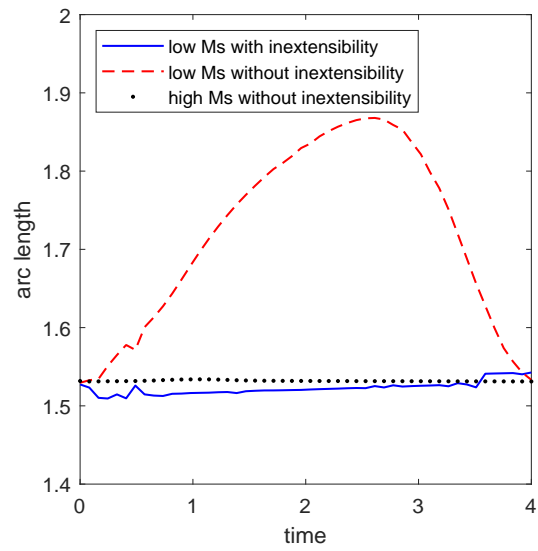}
		\caption{ Total arc length of vesicle membrane with the local inextensibility (blue line), and   the total arc lengths of vesicle membranes with low (100) (red dashed line) and high (20000) (black point) $\mathcal{M}_s$ and no local inextensibility, respectively, during vesicles passing through the constriction of the micro channel with otherwise same parameter values and settings.  }
		\label{fig:without lambda 2}
	\end{figure}

Effect of the local inextensibility of vesicle membrane is assessed by  comparing vesicle simulations with and without using the local inextensibility constraint $\mathcal {P}:\nabla \vel = 0$ in the model.  Snapshots of these simulations at different times are shown in Figure \ref{fig:contour}. They illustrate that a vesicle modeled without using the local inextensiblity can pass through the channel by introducing large extension and deformation of its body with a relatively small value of global inextensibility coefficient $\mathcal{M}_s$;  while a vesicle modeled with the local inextensibility hardly exhibits large extension and deformation of its body and blocks the channel.  This is also confirmed by Figure \ref{fig:without lambda 2}. It  shows under otherwise same conditions, the total arc length of the membrane of the vesicle modeled without the local inextensibility increases significantly when it passes through the channel, and the vesicle with the local inextensibilty preserves its membrane arc length well during the course of the simulation. 

\begin{figure}[h]
\centering     
\subfigure[] { \includegraphics[width=2.5in]{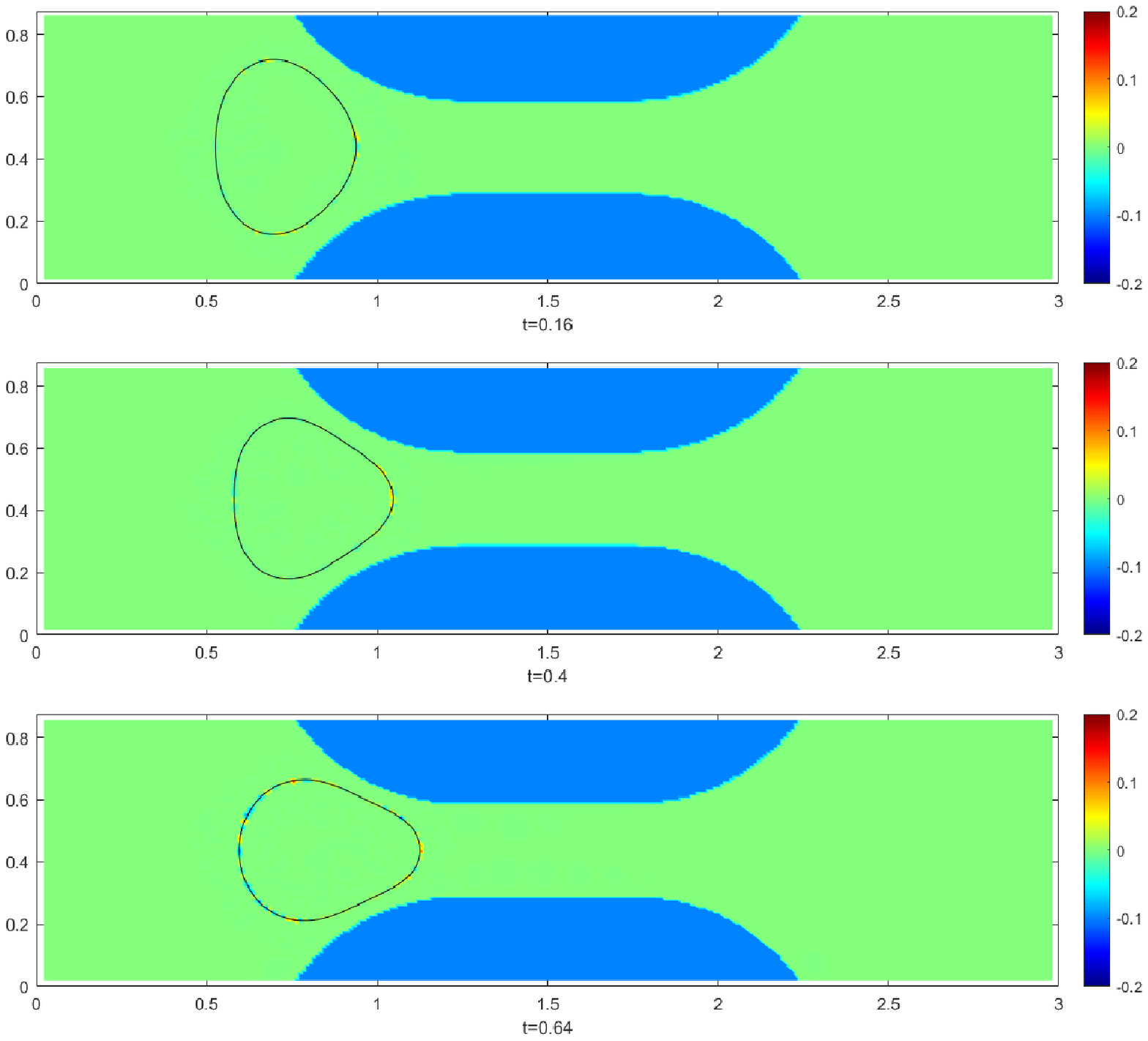}}
\subfigure [] {\includegraphics[width=2.5in]{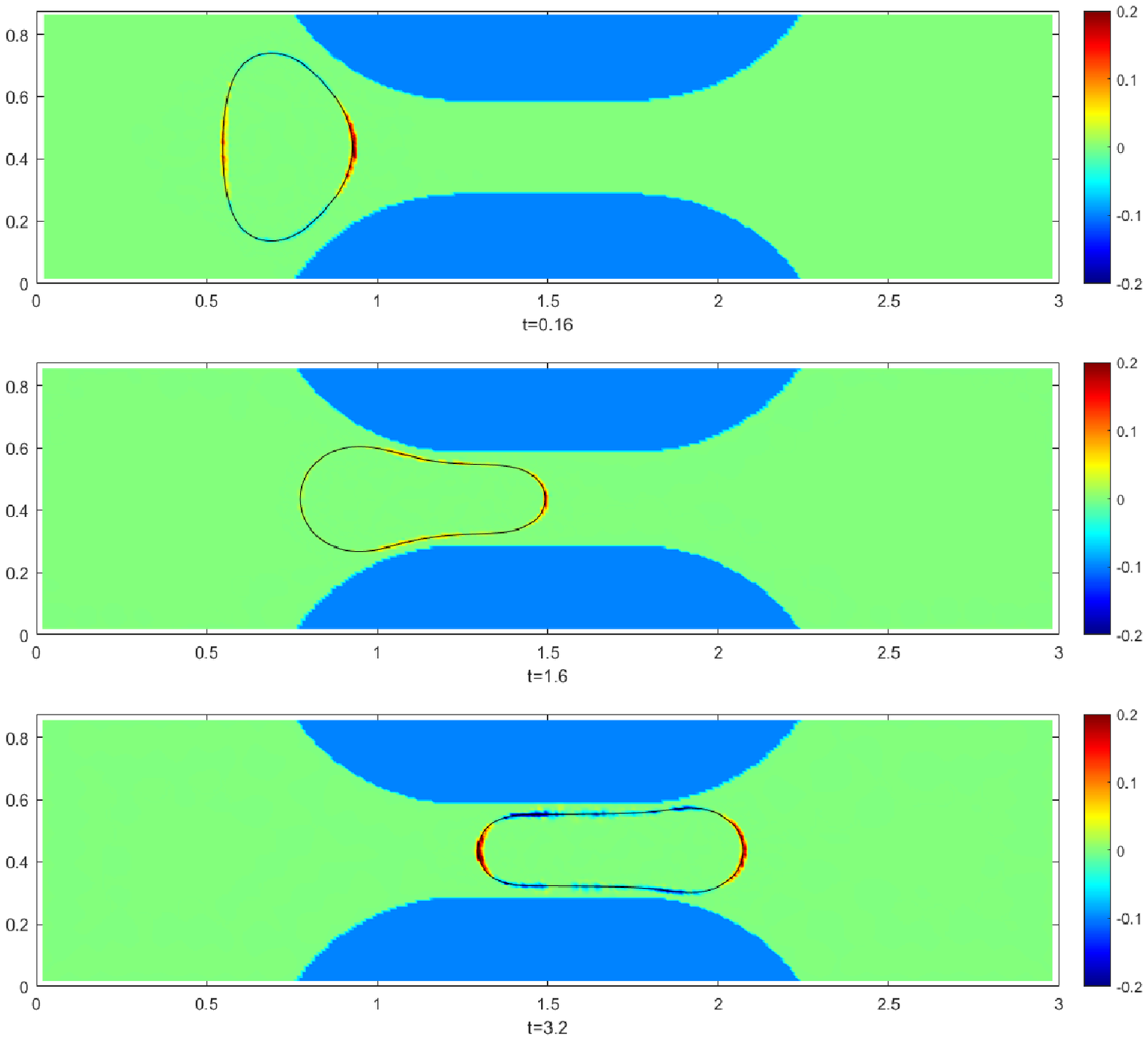}}
\subfigure[] { \includegraphics[width=2.5in]{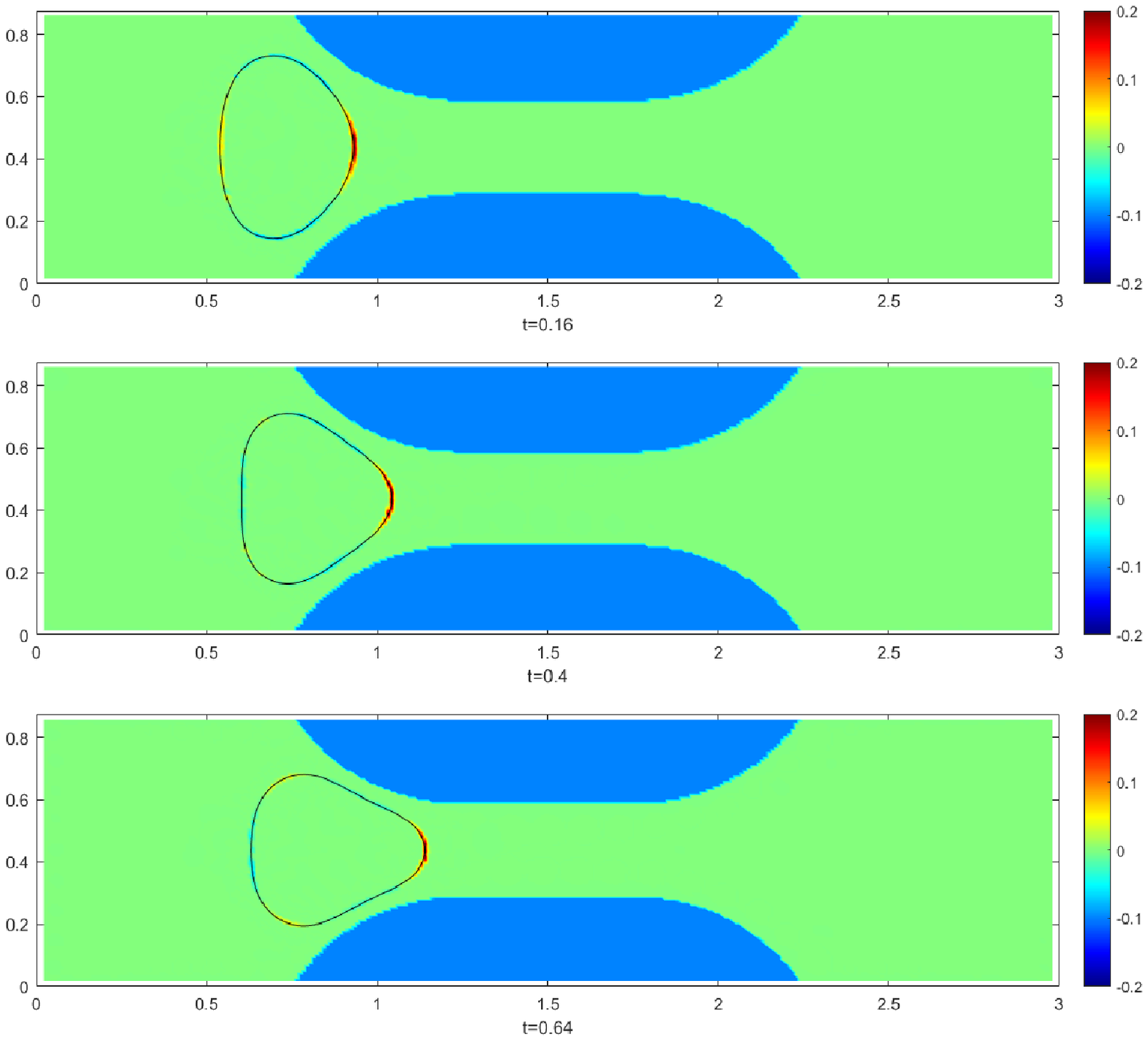}}
\caption{\label{fig:surface divergence}Effects of the local inextensibility $\mathcal{P}:\nabla\vel  =0$. Snapshots of membrane forces of vesicles: (a) $\mathcal{M}_s=100$ with the local inextensibility;  (b) $\mathcal{M}_s=100$ without the local inextensibility; and (c) $\mathcal{M}_s=20000$ without the local inextensibility. }
\end{figure}

Although the total arc length of a vesicle without the local inextensibility and with a very large $\mathcal{M}_s$ value could maintain almost unchanged  as shown in \ref{fig:contour} (c) and \ref{fig:without lambda 2}, the morphological changes of vesicles with and without the local inextensibility are drastically  different. 
For the vesicles modeled without the local inextensibility, 
 Figure \ref{fig:surface divergence}(b,c) illustrates that  the vesicle membranes are  stretched (red) or compressed (blue) everywhere, even though  the total arc length of the vesicle modeled using a large  modulus $\mathcal{M}_s$ value could be preserved, and the vesicle forms the blockage.  For the vesicle modeled with the local inextensibility,  Figure \ref{fig:surface divergence}(c)  confirms that there is almost no local extension or compression of the membrane, which is consistent with experimental observations. All simulations described below use the local inextensibility.

\begin{figure}[h]
		\centering
	\includegraphics[width=4.in]{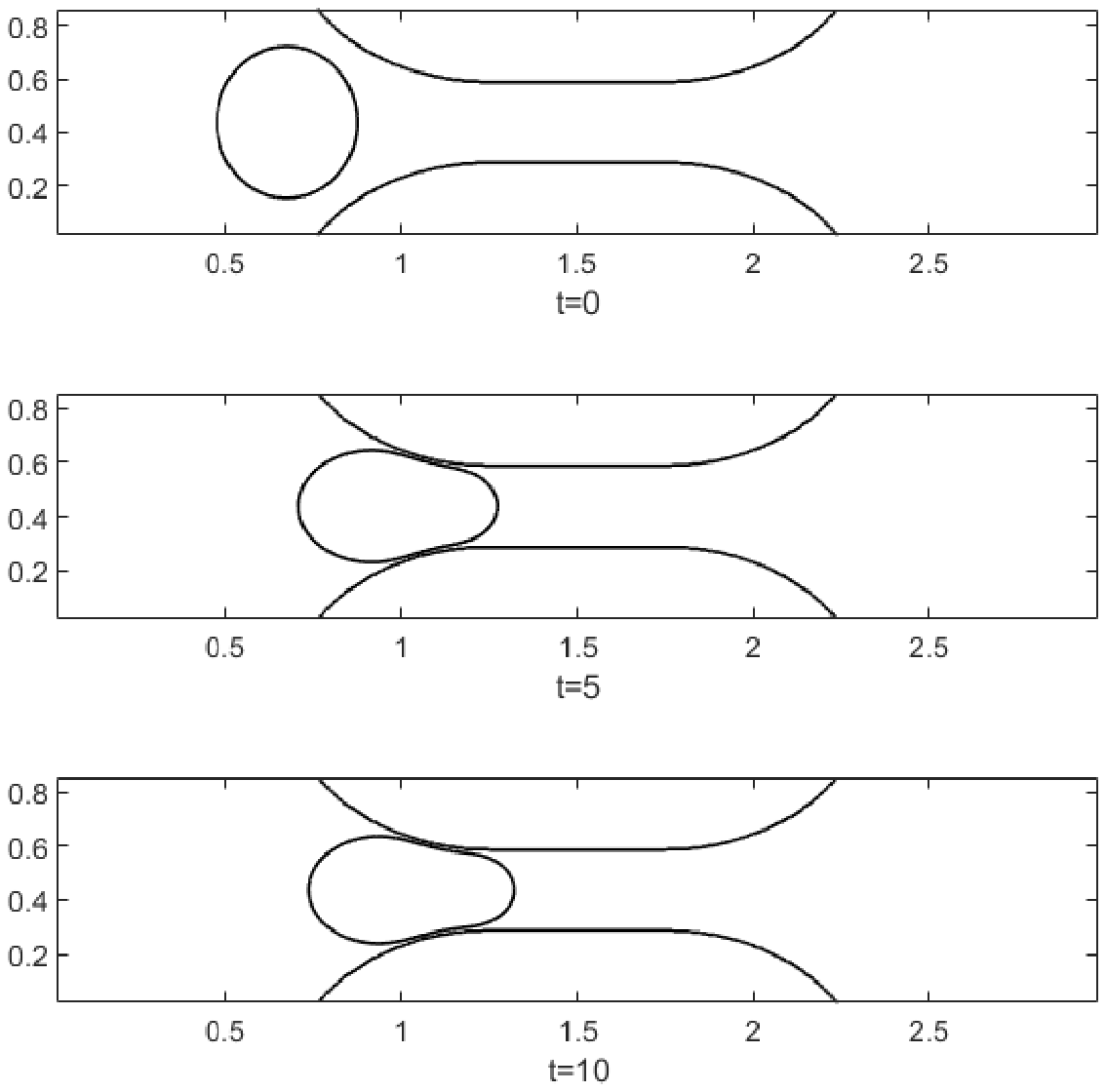}
		\caption{Side view of a vesicle with  surface-volume ratio $1.5:1$ at different times.}
		\label{fig:small ratio}
	\end{figure}

\begin{figure}[h]
		\centering
		\includegraphics[width=4.in]{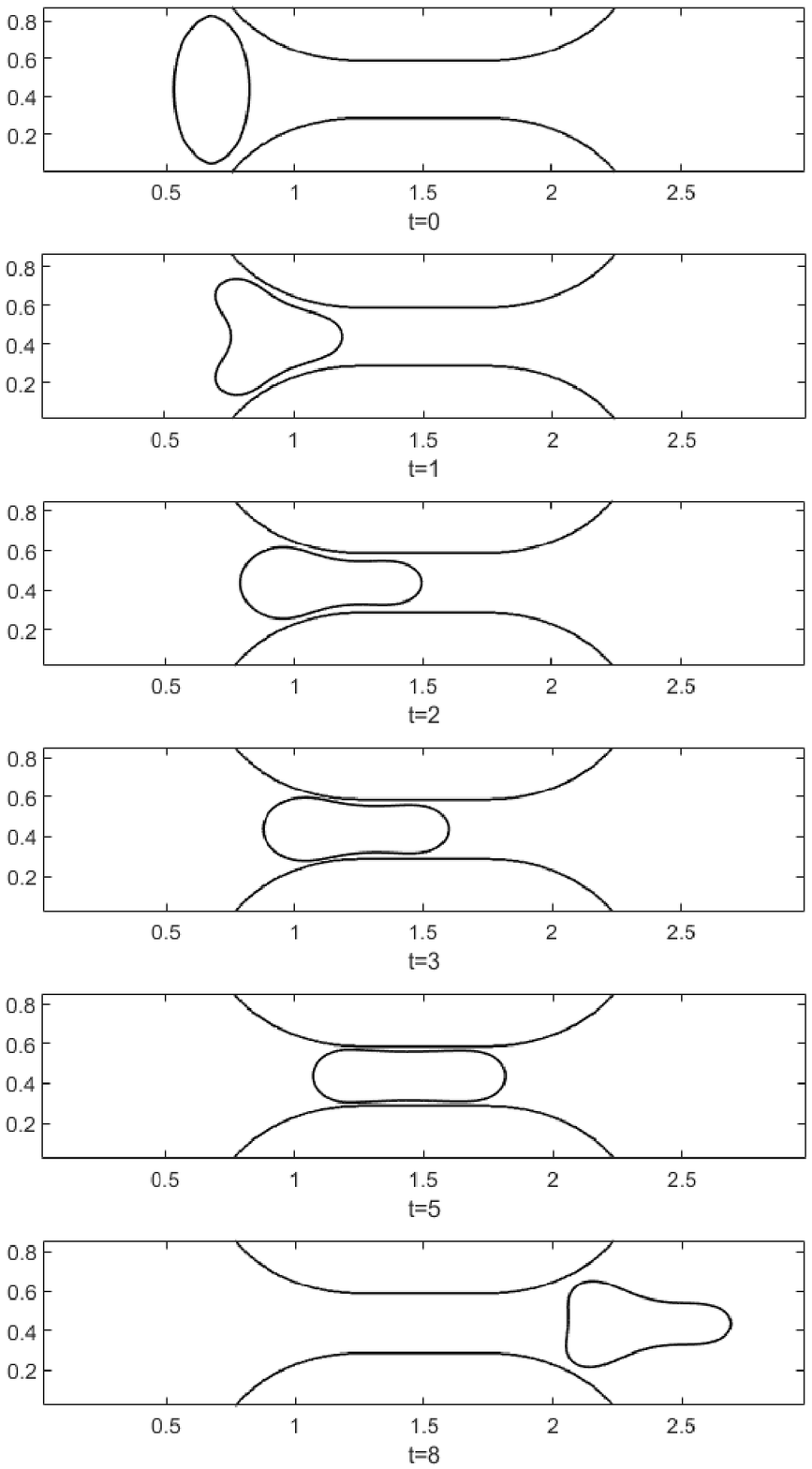} 
		\caption{Side view of a vesicle with surface-volume ratio $2:1$ at different times.}
		\label{fig:large ratio}
	\end{figure}

Both experiments and clinic reports have shown that  the cell bending modulus and surface-volume  ratio play important roles in determining the deformability of vesicles, especially when they pass through narrow channels \cite{takagi2009deformation,namvar2020surface,renoux2019impact}. The latest results reveal that a moderate decrease in  the surface-volume ratio  has a more significant effect than varying the  cell bending  stiffness. This surface-volume ratio effect is tested by increasing the ratio value slightly from  $ 1.5:1$ to $2:1$. Results in Figures \ref{fig:small ratio} and \ref{fig:large ratio} confirm that the more rounded vesicles are much harder to pass through the narrow channel and can easily form the blockage. This is consistent with the experimental observations.

 The effect of the bending modulus is assessed by increasing its value 10 times. The surface-volume ratio of the vesicle is $2:1$ in this test. 
 Figure \ref{fig:large ratio} illustrates that this more rigid vesicle  can also pass the same size channel but exhibits very different shape transformation.
	
\begin{figure}[h]
		\centering
		\includegraphics[width=4.in]{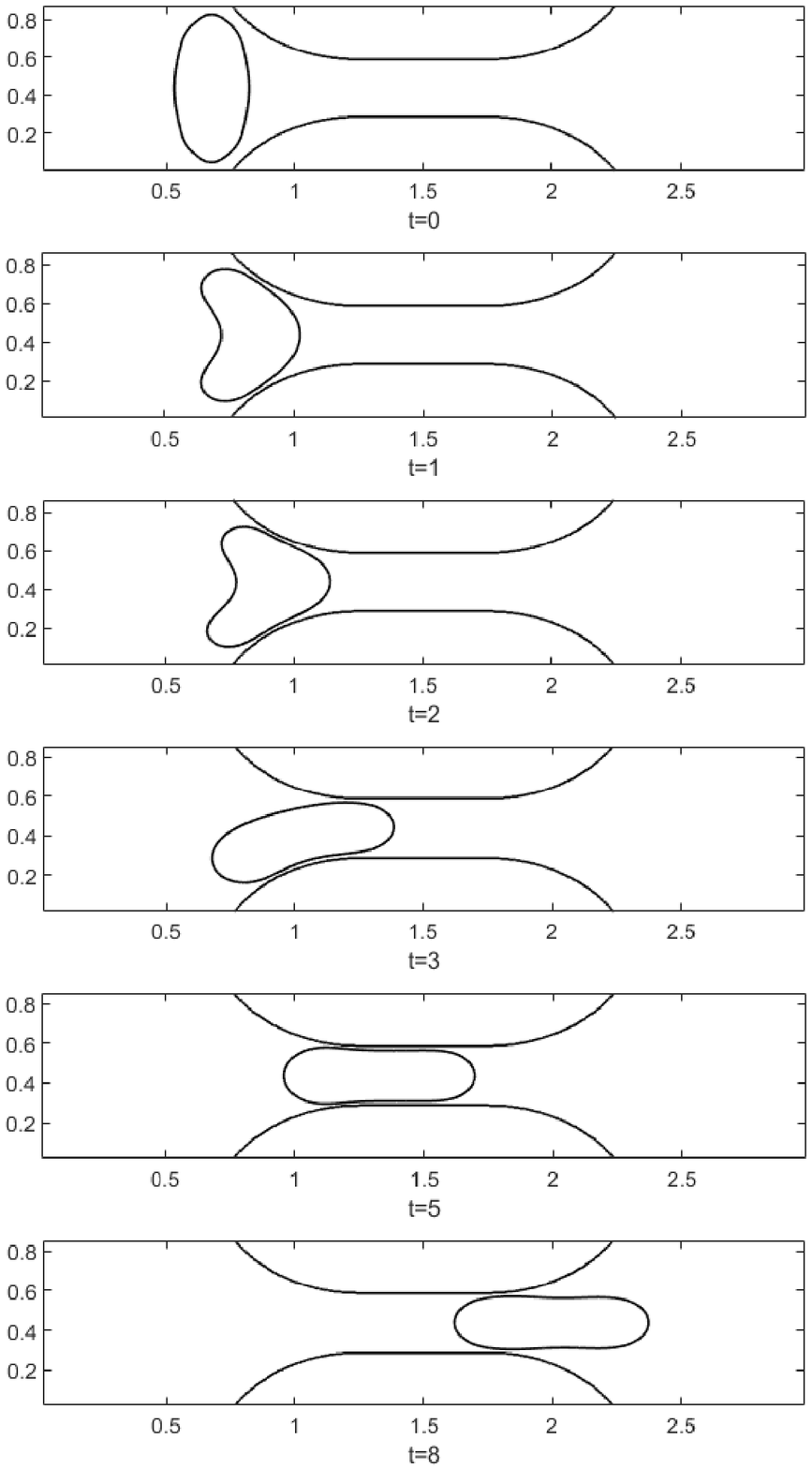} 
		\caption{Side view of a vesicle with  large bending modulus $\kappa_B=4\times 10^{-1}$ and surface-volume ratio $2:1$ at different times.}
		\label{fig:larger bending}
	\end{figure}

\section{Conclusion}

In this paper, an energy variational method is used to derive a thermodynamically consistent phase-field model for simulating vesicles motion and deformation. Corresponding Allen-Cahn GNBC boundary conditions accounting for vesicle-wall (or fluid-structure) interaction are also proposed by introducing the boundary dissipation and the vesicle-wall interaction energy.  

Then an efficient $C^0$ finite element method combined with the mid-point temporal discretization is proposed to solve the obtained model equations.  Thanks to the mid-point temporal discretization, the obtained numerical scheme is unconditionally energy stable. The numerical experiments confirm that this scheme is second-order accurate in both space and time. The vesicle tank treading and tumbling simulations reproduce experimental observations. And the flipping ellipse simulation agrees with the analytical solution well. Finally, the model is used to investigate how  vesicles' mechanical properties affect the vesicles' capability to pass through narrow channels. It is shown that whether a vesicle can pass a narrow channel is largely determined by  the surface-volume ratio of the vesicle, which is consistent with {\it in vitro} experiments.

\section*{Authors’ contributions}

L.S, Z.X, and S.X did the model derivations. L.S, P.L and S.X carried out the numerical analysis and simulations. Z.X, H.H, P.L and S.X designed  and coordinated the study. All participated in the preparation of the manuscript. All authors gave final approval for publication.

\section*{Declaration of competing interest}
The authors declare that they have no known competing financial interests or personal relationships that could have appeared to influence the work reported in this paper.

\section*{Acknowledgements}
This work was partially supported by the NSFC (grant numbers 12071190, 11771040, 11861131004, 91430106), and NSERC (CA). L.S was partially supported by the Chinese Scholarship Council for studying at the University of Dundee. Z.X was partially supported by
the NSF CDS\&E-MSS 1854779 and NSF-1821242.

\newpage

\newpage
\bibliographystyle{elsarticle-num} 
\bibliography{Mybib}

\begin{thebibliography}{10}
\expandafter\ifx\csname url\endcsname\relax
  \def\url#1{\texttt{#1}}\fi
\expandafter\ifx\csname urlprefix\endcsname\relax\def\urlprefix{URL }\fi
\expandafter\ifx\csname href\endcsname\relax
  \def\href#1#2{#2} \def\path#1{#1}\fi

\bibitem{Takeishi_Ito_2019}
N.~Takeishi, H.~Ito, M.~Kaneko, S.~Wada, Deformation of a red blood cell in a
  narrow rectangular microchannel, Micromachines 10~(3) (2019) 199.

\bibitem{Caimi_Presti_04}
C.~G, P.~RL., Techniques to evaluate erythrocyte deformability in diabetes
  mellitus, Acta Diabetol. 41~(3) (2004) 99--103.

\bibitem{Pschl2003EndotoxinBT}
J.~P{\"o}schl, C.~Leray, P.~Ruef, J.~Cazenave, O.~Linderkamp, Endotoxin binding
  to erythrocyte membrane and erythrocyte deformability in human sepsis and in
  vitro, Critical Care Medicine 31 (2003) 924--928.

\bibitem{shin2017platelet}
E.-K. Shin, H.~Park, J.-Y. Noh, K.-M. Lim, J.-H. Chung, Platelet shape changes
  and cytoskeleton dynamics as novel therapeutic targets for anti-thrombotic
  drugs, Biomolecules \& therapeutics 25~(3) (2017) 223.

\bibitem{aslan2012platelet}
J.~E. Aslan, A.~Itakura, J.~M. Gertz, O.~J. McCarty, Platelet shape change and
  spreading, in: Platelets and Megakaryocytes, Springer, 2012, pp. 91--100.

\bibitem{lishuwang_2017_vesicle}
K.~Liu, G.~R. Marple, J.~Allard, S.~Li, S.~Veerapaneni, J.~Lowengrub, Dynamics
  of a multicomponent vesicle in shear flow, Soft matter 13~(19) (2017)
  3521--3531.

\bibitem{yangxiaofeng_2015_decoupled}
R.~Chen, G.~Ji, X.~Yang, H.~Zhang, Decoupled energy stable schemes for
  phase-field vesicle membrane model, Journal of Computational Physics 302
  (2015) 509--523.

\bibitem{duqiang_2004_bending}
Q.~Du, C.~Liu, X.~Wang, A phase field approach in the numerical study of the
  elastic bending energy for vesicle membranes, Journal of Computational
  Physics 198~(2) (2004) 450--468.

\bibitem{Noguchi14159}
H.~Noguchi, G.~Gompper, \href{https://www.pnas.org/content/102/40/14159}{Shape
  transitions of fluid vesicles and red blood cells in capillary flows},
  Proceedings of the National Academy of Sciences 102~(40) (2005) 14159--14164.
\newblock \href
  {http://arxiv.org/abs/https://www.pnas.org/content/102/40/14159.full.pdf}
  {\path{arXiv:https://www.pnas.org/content/102/40/14159.full.pdf}}, \href
  {https://doi.org/10.1073/pnas.0504243102}
  {\path{doi:10.1073/pnas.0504243102}}.
\newline\urlprefix\url{https://www.pnas.org/content/102/40/14159}

\bibitem{li_2018_dpd}
Z.~Li, X.~Bian, Y.-H. Tang, G.~E. Karniadakis, A dissipative particle dynamics
  method for arbitrarily complex geometries, Journal of Computational Physics
  355 (2018) 534--547.

\bibitem{pivkin2009effect}
I.~V. Pivkin, P.~D. Richardson, G.~E. Karniadakis, Effect of red blood cells on
  platelet aggregation, IEEE Engineering in Medicine and Biology Magazine
  28~(2) (2009) 32--37.

\bibitem{pivkin2008accurate}
I.~V. Pivkin, G.~E. Karniadakis, Accurate coarse-grained modeling of red blood
  cells, Physical review letters 101~(11) (2008) 118105.

\bibitem{xu_level-set_2006}
J.-J. Xu, Z.~Li, J.~Lowengrub, H.~Zhao, A level-set method for interfacial
  flows with surfactant, Journal of Computational Physics 212~(2) (2006)
  590--616.
\newblock \href {https://doi.org/10.1016/j.jcp.2005.07.016}
  {\path{doi:10.1016/j.jcp.2005.07.016}}.

\bibitem{xu_level-set_2014}
J.-J. Xu, W.~Ren, A level-set method for two-phase flows with moving contact
  line and insoluble surfactant, Journal of Computational Physics 263 (2014)
  71--90.
\newblock \href {https://doi.org/10.1016/j.jcp.2014.01.012}
  {\path{doi:10.1016/j.jcp.2014.01.012}}.

\bibitem{salac2011level}
D.~Salac, M.~Miksis, A level set projection model of lipid vesicles in general
  flows, Journal of Computational Physics 230~(22) (2011) 8192--8215.

\bibitem{bonito2010parametric}
A.~Bonito, R.~H. Nochetto, M.~S. Pauletti, Parametric fem for geometric
  biomembranes, Journal of Computational Physics 229~(9) (2010) 3171--3188.

\bibitem{jenkins1977equations}
J.~T. Jenkins, The equations of mechanical equilibrium of a model membrane,
  SIAM Journal on Applied Mathematics 32~(4) (1977) 755--764.

\bibitem{hu2007continuum}
D.~Hu, P.~Zhang, E.~Weinan, Continuum theory of a moving membrane, Physical
  Review E 75~(4) (2007) 041605.

\bibitem{hao2015fictitious}
W.~Hao, Z.~Xu, C.~Liu, G.~Lin, A fictitious domain method with a hybrid cell
  model for simulating motion of cells in fluid flow, Journal of Computational
  Physics 280 (2015) 345--362.

\bibitem{laimingzhi_2010_immersed_boundary}
Y.~Kim, M.-C. Lai, Simulating the dynamics of inextensible vesicles by the
  penalty immersed boundary method, Journal of Computational Physics 229~(12)
  (2010) 4840--4853.

\bibitem{laimingzhi_2016_immersed_boundary}
Y.~Seol, W.-F. Hu, Y.~Kim, M.-C. Lai, An immersed boundary method for
  simulating vesicle dynamics in three dimensions, Journal of Computational
  Physics 322 (2016) 125--141.

\bibitem{laimingzhi_2020_immersed_boundary}
K.~C. Ong, M.-C. Lai, An immersed boundary projection method for simulating the
  inextensible vesicle dynamics, Journal of Computational Physics 408 (2020)
  109277.

\bibitem{wang2020immersed}
X.~Wang, X.~Gong, K.~Sugiyama, S.~Takagi, H.~Huang, An immersed boundary method
  for mass transfer through porous biomembranes under large deformations,
  Journal of Computational Physics (2020) 109444.

\bibitem{wu2013simulation}
T.~Wu, J.~J. Feng, Simulation of malaria-infected red blood cells in
  microfluidic channels: Passage and blockage, Biomicrofluidics 7~(4) (2013)
  044115.

\bibitem{hu2016vesicle}
W.-F. Hu, M.-C. Lai, Y.~Seol, Y.-N. Young, Vesicle electrohydrodynamic
  simulations by coupling immersed boundary and immersed interface method,
  Journal of Computational Physics 317 (2016) 66--81.

\bibitem{kolahdouz2015numerical}
E.~M. Kolahdouz, D.~Salac, A numerical model for the trans-membrane voltage of
  vesicles, Applied Mathematics Letters 39 (2015) 7--12.

\bibitem{duqiang_2009_variational}
Q.~Du, C.~Liu, R.~Ryham, X.~Wang, Energetic variational approaches in modeling
  vesicle and fluid interactions, Physica D: Nonlinear Phenomena 238~(9-10)
  (2009) 923--930.

\bibitem{lin2006simulations}
P.~Lin, C.~Liu, Simulations of singularity dynamics in liquid crystal flows: A
  c0 finite element approach, Journal of Computational Physics 215~(1) (2006)
  348--362.

\bibitem{Diegel2019finite}
A.~E. Diegel, S.~W. Walker, A finite element method for a phase field model of
  nematic liquid crystal droplets, Communications in Computational Physics 25
  (2019) 155--188.

\bibitem{wise2009energy}
S.~M. Wise, C.~Wang, J.~S. Lowengrub, An energy-stable and convergent
  finite-difference scheme for the phase field crystal equation, SIAM Journal
  on Numerical Analysis 47~(3) (2009) 2269--2288.

\bibitem{xu2018three}
S.~Xu, M.~Alber, Z.~Xu, Three-phase model of visco-elastic incompressible fluid
  flow and its computational implementation, Communications in Computational
  Physics. 25(2) (2019) 586--624.

\bibitem{wangqi_2011_finite_element}
J.~Hua, P.~Lin, C.~Liu, Q.~Wang, Energy law preserving c0 finite element
  schemes for phase field models in two-phase flow computations, Journal of
  Computational Physics 230~(19) (2011) 7115--7131.

\bibitem{yangxiaofeng_2017_3_component}
X.~Yang, J.~Zhao, Q.~Wang, J.~Shen, Numerical approximations for a
  three-component cahn--hilliard phase-field model based on the invariant
  energy quadratization method, Mathematical Models and Methods in Applied
  Sciences 27~(11) (2017) 1993--2030.

\bibitem{lowengrub2009phase}
J.~S. Lowengrub, A.~R{\"a}tz, A.~Voigt, Phase-field modeling of the dynamics of
  multicomponent vesicles: Spinodal decomposition, coarsening, budding, and
  fission, Physical Review E 79~(3) (2009) 031926.

\bibitem{du2005phase}
Q.~Du, C.~Liu, R.~Ryham, X.~Wang, A phase field formulation of the willmore
  problem, Nonlinearity 18~(3) (2005) 1249.

\bibitem{du2007analysis}
Q.~Du, M.~Li, C.~Liu, Analysis of a phase field navier-stokes vesicle-fluid
  interaction model, Discrete \& Continuous Dynamical Systems-B 8~(3) (2007)
  539.

\bibitem{duqiang_2005_curvature}
Q.~Du, C.~Liu, R.~Ryham, X.~Wang, Modeling the spontaneous curvature effects in
  static cell membrane deformations by a phase field formulation,
  Communications on Pure \& Applied Analysis 4~(3) (2005) 537.

\bibitem{biben2005phase}
T.~Biben, K.~Kassner, C.~Misbah, Phase-field approach to three-dimensional
  vesicle dynamics, Physical Review E 72~(4) (2005) 041921.

\bibitem{zhang2009phase}
J.~Zhang, S.~Das, Q.~Du, A phase field model for vesicle--substrate adhesion,
  Journal of Computational Physics 228~(20) (2009) 7837--7849.

\bibitem{chen2015decoupled}
R.~Chen, G.~Ji, X.~Yang, H.~Zhang, Decoupled energy stable schemes for
  phase-field vesicle membrane model, Journal of Computational Physics 302
  (2015) 509--523.

\bibitem{wu2013strong}
H.~Wu, X.~Xu, Strong solutions, global regularity, and stability of a
  hydrodynamic system modeling vesicle and fluid interactions, SIAM Journal on
  Mathematical Analysis 45~(1) (2013) 181--214.

\bibitem{duqiang_2008_bending}
Q.~Du, J.~Zhang, Adaptive finite element method for a phase field bending
  elasticity model of vesicle membrane deformations, SIAM Journal on Scientific
  Computing 30~(3) (2008) 1634--1657.

\bibitem{voigt_2014_local_inextensibility}
S.~Aland, S.~Egerer, J.~Lowengrub, A.~Voigt, Diffuse interface models of
  locally inextensible vesicles in a viscous fluid, Journal of computational
  physics 277 (2014) 32--47.

\bibitem{beaucourt_2004_steady}
J.~Beaucourt, F.~Rioual, T.~S{\'e}on, T.~Biben, C.~Misbah, Steady to unsteady
  dynamics of a vesicle in a flow, Physical Review E 69~(1) (2004) 011906.

\bibitem{shen2020energy}


\bibitem{xu2014energetic}
S.~Xu, P.~Sheng, C.~Liu, An energetic variational approach for ion transport,
  Communications in Mathematical Sciences 12~(4) (2014) 779--789.

\bibitem{xu2018osmosis}
S.~Xu, B.~Eisenberg, Z.~Song, H.~Huang, Osmosis through a semi-permeable
  membrane: a consistent approach to interactions, arXiv preprint
  arXiv:1806.00646 (2018).

\bibitem{guo2015thermodynamically}
Z.~Guo, P.~Lin, A thermodynamically consistent phase-field model for two-phase
  flows with thermocapillary effects, Journal of Fluid Mechanics 766 (2015)
  226--271.

\bibitem{qian_2006_slipBC}
T.~Qian, X.-P. Wang, P.~Sheng, A variational approach to the moving contact
  line hydrodynamics, arXiv preprint cond-mat/0602293 (2006).

\bibitem{cheng2019energy}
K.~Cheng, W.~Feng, C.~Wang, S.~M. Wise, An energy stable fourth order finite
  difference scheme for the cahn--hilliard equation, Journal of Computational
  and Applied Mathematics 362 (2019) 574--595.

\bibitem{yan2018second}
Y.~Yan, W.~Chen, C.~Wang, S.~M. Wise, A second-order energy stable bdf
  numerical scheme for the cahn-hilliard equation, Commun. Comput. Phys. 23~(2)
  (2018) 572--602.

\bibitem{chen2016convergence}
W.~Chen, Y.~Liu, C.~Wang, S.~Wise, Convergence analysis of a fully discrete
  finite difference scheme for the cahn-hilliard-hele-shaw equation,
  Mathematics of Computation 85~(301) (2016) 2231--2257.

\bibitem{guo_2014_midpoint}
Z.~Guo, P.~Lin, J.~S. Lowengrub, A numerical method for the
  quasi-incompressible cahn--hilliard--navier--stokes equations for variable
  density flows with a discrete energy law, Journal of Computational Physics
  276 (2014) 486--507.

\bibitem{yang2017numerical}
X.~Yang, J.~Zhao, Q.~Wang, Numerical approximations for the molecular beam
  epitaxial growth model based on the invariant energy quadratization method,
  Journal of Computational Physics 333 (2017) 104--127.

\bibitem{yang2020convergence}
X.~Yang, G.-D. Zhang, Convergence analysis for the invariant energy
  quadratization (ieq) schemes for solving the cahn--hilliard and allen--cahn
  equations with general nonlinear potential, Journal of Scientific Computing
  82~(3) (2020) 1--28.

\bibitem{shen2018convergence}
J.~Shen, J.~Xu, Convergence and error analysis for the scalar auxiliary
  variable (sav) schemes to gradient flows, SIAM Journal on Numerical Analysis
  56~(5) (2018) 2895--2912.

\bibitem{shen2018scalar}
J.~Shen, J.~Xu, J.~Yang, The scalar auxiliary variable (sav) approach for
  gradient flows, Journal of Computational Physics 353 (2018) 407--416.

\bibitem{shen2019new}
J.~Shen, J.~Xu, J.~Yang, A new class of efficient and robust energy stable
  schemes for gradient flows, SIAM Review 61~(3) (2019) 474--506.

\bibitem{gong2019arbitrarily}
Y.~Gong, J.~Zhao, Q.~Wang, Arbitrarily high-order unconditionally energy stable
  sav schemes for gradient flow models, Computer Physics Communications (2019)
  107033.

\bibitem{li2020sav}
X.~Li, J.~Shen, On a sav-mac scheme for the cahn-hilliard-navier-stokes phase
  field model and its error analysis for the corresponding cahn-hilliard-stokes
  case, Mathematical Models and Methods in Applied Sciences (2020).

\bibitem{gu2016two}
R.~Gu, X.~Wang, M.~Gunzburger, A two phase field model for tracking
  vesicle--vesicle adhesion, Journal of mathematical biology 73~(5) (2016)
  1293--1319.

\bibitem{gu2014simulating}
R.~Gu, X.~Wang, M.~Gunzburger, Simulating vesicle--substrate adhesion using two
  phase field functions, Journal of Computational Physics 275 (2014) 626--641.

\bibitem{guillen2018unconditionally}
F.~Guill{\'e}n-Gonz{\'a}lez, G.~Tierra, Unconditionally energy stable numerical
  schemes for phase-field vesicle membrane model, Journal of computational
  physics 354 (2018) 67--85.

\bibitem{marth2016margination}
W.~Marth, S.~Aland, A.~Voigt, Margination of white blood cells: a computational
  approach by a hydrodynamic phase field model, Journal of Fluid Mechanics 790
  (2016) 389--406.

\bibitem{qian_variational_2006}
T.~Qian, X.-P. Wang, P.~Sheng, A variational approach to moving contact line
  hydrodynamics, Journal of Fluid Mechanics 564 (2006) 333.
\newblock \href {https://doi.org/10.1017/S0022112006001935}
  {\path{doi:10.1017/S0022112006001935}}.

\bibitem{eisenberg2010energy}
B.~Eisenberg, Y.~Hyon, C.~Liu, Energy variational analysis of ions in water and
  channels: Field theory for primitive models of complex ionic fluids, The
  Journal of Chemical Physics 133~(10) (2010) 104104.

\bibitem{lee2012regularized}
H.~G. Lee, J.~Kim, Regularized dirac delta functions for phase field models,
  International journal for numerical methods in engineering 91~(3) (2012)
  269--288.

\bibitem{qian_molecular_2006}
T.~Qian, X.-P. Wang, P.~Sheng, Molecular hydrodynamics of the moving contact
  line in two-phase immersible flows, Commun. Comput. Phys. (2006) 52.

\bibitem{ren2007boundary}
W.~Ren, W.~E, Boundary conditions for the moving contact line problem, Physics
  of fluids 19~(2) (2007) 022101.

\bibitem{ren_heterogeneous_2005}
W.~Ren, W.~E, Heterogeneous multiscale method for the modeling of complex
  fluids and micro-fluidics, Journal of Computational Physics 204~(1) (2005)
  1--26.
\newblock \href {https://doi.org/10.1016/j.jcp.2004.10.001}
  {\path{doi:10.1016/j.jcp.2004.10.001}}.

\bibitem{laadhari2012vesicle}
A.~Laadhari, P.~Saramito, C.~Misbah, Vesicle tumbling inhibited by inertia,
  Physics of Fluids 24~(3) (2012) 031901.

\bibitem{basu2011tank}
H.~Basu, A.~K. Dharmadhikari, J.~A. Dharmadhikari, S.~Sharma, D.~Mathur, Tank
  treading of optically trapped red blood cells in shear flow, Biophysical
  journal 101~(7) (2011) 1604--1612.

\bibitem{fischer2004shape}
T.~M. Fischer, Shape memory of human red blood cells, Biophysical journal
  86~(5) (2004) 3304--3313.

\bibitem{jeffery1922motion}
G.~B. Jeffery, The motion of ellipsoidal particles immersed in a viscous fluid,
  Proceedings of the Royal Society of London. Series A, Containing papers of a
  mathematical and physical character 102~(715) (1922) 161--179.

\bibitem{han2019flow}
Y.~Han, H.~Lin, M.~Ding, R.~Li, T.~Shi, Flow-induced translocation of vesicles
  through a narrow pore, Soft matter 15~(16) (2019) 3307--3314.

\bibitem{takagi2009deformation}
S.~Takagi, T.~Yamada, X.~Gong, Y.~Matsumoto, The deformation of a vesicle in a
  linear shear flow, Journal of applied mechanics 76~(2) (2009).

\bibitem{namvar2020surface}
A.~Namvar, A.~J. Blanch, M.~W. Dixon, O.~M. Carmo, B.~Liu, S.~Tiash, O.~Looker,
  D.~Andrew, L.-J. Chan, W.-H. Tham, et~al., Surface area-to-volume ratio, not
  cellular viscoelasticity, is the major determinant of red blood cell
  traversal through small channels, Cellular Microbiology (2020) e13270.

\bibitem{renoux2019impact}
C.~Renoux, M.~Faivre, A.~Bessaa, L.~Da~Costa, P.~Joly, A.~Gauthier, P.~Connes,
  Impact of surface-area-to-volume ratio, internal viscosity and membrane
  viscoelasticity on red blood cell deformability measured in isotonic
  condition, Scientific reports 9~(1) (2019) 1--7.

\end{thebibliography}

\end{document}